\documentclass{amsart}

\usepackage[colorlinks=true, urlcolor=black, citecolor=black, linkcolor=black, hyperfootnotes=true]{hyperref}
\usepackage{amssymb}
\usepackage{mathrsfs}
\usepackage{aliascnt}

\numberwithin{equation}{section}

\newtheorem{thm}{Theorem}[section]

\newaliascnt{prp}{thm}
\newtheorem{prp}[prp]{Proposition}
\aliascntresetthe{prp}

\newtheorem*{prp*}{Proposition}

\newaliascnt{cor}{thm}
\newtheorem{cor}[cor]{Corollary}
\aliascntresetthe{cor}

\theoremstyle{definition}

\newaliascnt{dfn}{thm}
\newtheorem{dfn}[dfn]{Definition}
\aliascntresetthe{dfn}

\newaliascnt{xpl}{thm}
\newtheorem{xpl}[xpl]{Example}
\aliascntresetthe{xpl}

\newaliascnt{rmk}{thm}
\newtheorem{rmk}[rmk]{Remark}
\aliascntresetthe{rmk}

\author{Tristan Bice}
\thanks{The author is supported by the GA\v{C}R project EXPRO 20-31529X and RVO: 67985840.}
\address{Institute of Mathematics of the Czech Academy of Sciences, \v{Z}itn\'a 25, Prague}
\email{bice@math.cas.cz}

\title[Representing Structured Semigroups on \'Etale Groupoid Bundles]{Representing Structured Semigroups\\ on \'Etale Groupoid Bundles}

\begin{document}

\begin{abstract}
We examine a semigroup analogue of the Kumjian-Renault representation of C*-algebras with Cartan subalgebras on twisted groupoids.  Specifically, we show how to represent semigroups with distinguished normal subsemigroups as `slice-sections' of groupoid bundles.
\end{abstract}

\maketitle

\section*{Introduction}

\subsection*{Background}

A cornerstone of groupoid C*-algebra theory is the Kumjian-Renault Weyl groupoid representation from \cite{Kumjian1986} and \cite{Renault2008}.  Specifically, given a C*-algebra $A$ with a Cartan subalgebra $C$, they showed how to construct a Fell line bundle over an \'etale groupoid on which to represent $A$ as continuous sections, thus realising the abstract algebra $A$ as a concrete twisted groupoid C*-algebra.

One key aspect of their construction is the use of normalisers of $C$ to construct an appropriate groupoid of germs.  Recently it has become apparent that this normaliser semigroup $S$ already contains a significant amount of structural information.  For example, if $A$ is already known to be groupoid C*-algebra then this groupoid can be recovered just from the semigroup structure of $S$ and $C$ (as the groupoid of ultrafilters w.r.t. the `domination relation' on $S$ \textendash\, see \cite{BiceClark2021}).

This begs the question of whether one can similarly represent a given abstract semigroup $S$, together with some additional structure coming from appropriate subsemigroups, as a semigroup of sections of some groupoid bundle.  This is borne out in the present paper, where we construct universal representations of $S$ on groupoid bundles formed from cosets and filters of $S$.  The basic idea is to extend the relevant theory of inverse semigroups, replacing the idempotents with normal subsemigroups analogous to Cartan subalgebras.

\subsection*{Outline}

As a gentle introduction to our more general results, in \autoref{Motivation} we first consider simpler classes of semigroups.  We start with semilattices and their filters in \autoref{Semilattices}, moving on to inverse semigroups and their cosets in \autoref{InverseSemigroups}.  Finally we introduce \emph{structured semigroups} in \autoref{StructuredSemigroups}, a generalisation of inverse semigroups that encompasses our primary motivating examples, namely the normaliser semigroups of Cartan subalgebras of C*-algebras.

In \autoref{SliceSections} we introduce groupoid bundles and examine their slice-sections.  We give some further motivational comments at the end to indicate how we are going to construct groupoid bundles on which to represent structured semigroups.

Next in \autoref{Morphisms} we discuss appropriate morphisms between \'etale groupoids and groupoid bundles.  Specifically, we consider relational Zakrzewski morphisms between \'etale groupoids, which are needed to include both continuous functions and group homomorphisms (but in opposite directions \textendash\, see \autoref{SpacesVsGroups}), as well as the filter-coset Zakrzewski morphism $\vartriangleleft$ appearing later in \autoref{Universality}.  Building on these, we consider Zakrzewski-Pierce morphisms between groupoid bundles, which we show in \autoref{BundleMorphism->SemigroupHomo} correspond to semigroup homomorphisms of their slice-sections.

From \autoref{Domination}, we make the standing assumption that $(S,N,Z)$ is a structured semigroup.  We then introduce the domination relation, a generalisation of the usual ordering on an inverse semigroup which also has analogous properties that we proceed to examine.  Next in \autoref{Duals} we introduce dual subsets as something of a substitute for inverses.  Optimal results here require an additional diagonality assumption on $N$ that we discuss in \autoref{Diagonality}, which will become more relevant again later in \autoref{Filters}.

Continuing on the path of extending inverse semigroup theory, we introduce atlases and cosets in \autoref{Atlases}.  In \autoref{Cosets}, we show that the cosets form an \'etale groupoid -- see \autoref{EtaleCosets}.  Again optimal results require an additional assumption, namely that $Z$ is `symmetric', as discussed in \autoref{Symmetry}.  Under this assumption, we use cosets to construct a universal \'etale representation in \autoref{SymmetricEtaleRepresentation}.

In \autoref{Equivalence} and \autoref{Bundles}, we show how to split up cosets into certain equivalence classes.  The main result is \autoref{CosetBundleEtale}, which says that these form an \'etale bundle.  Next we represent $S$ on this coset bundle in \autoref{TheUniversalRepresentation}.  The main result here is \autoref{CosetBundleUniversal}, which says that this is a universal bundle representation, assuming symmetry.

Finally in \autoref{Filters}, we take a closer look at filters or, more accurately, directed cosets.  These are shown to form an \'etale subgroupoid of the coset groupoid which instead yield Zakrzewski-universal representations (again under symmetry), as shown in \autoref{ZakrzewskiUniversalEtaleRepresentation} and \autoref{ZakrzewskiUniversalBundleRepresentation}.

\subsection*{Acknowledgements}

Many thanks to the anonymous referee for reading the first draft very carefully and making numerous helpful comments and suggestions which have been incorporated into this final version.

\section{Motivation}\label{Motivation}

\subsection{Semilattices}\label{Semilattices}

To put our results in context, it is instructive to first consider semilattices, i.e. commutative semigroups consisting entirely of idempotents.  For example, the open subsets $\mathcal{O}(X)$ of any topological space $X$ form a semilattice under the intersection operation $\cap$.  Of course the same applies to any subfamily of $\mathcal{O}(X)$ that is closed under pairwise intersections, like the family of open intervals $(a,b)$ on the real line $\mathbb{R}$.

A slightly less trivial observation is that, in fact, all semilattices are of this form, at least up to isomorphism.  Put another way, every semilattice has a faithful spatial representation.  Formally, a \emph{spatial representation} of a semilattice $S$ on a space $X$ is a map $\theta:S\rightarrow\mathcal{O}(X)$ such that $\theta[S]$ covers $X$ and, for all $a,b\in S$,
\[\theta(ab)=\theta(a)\cap\theta(b).\]
To call a representation \emph{faithful} just means it is injective, i.e. $\theta(a)=\theta(b)$ iff $a=b$.

To construct a faithful representation of any given semilattice $S$, we can consider its \emph{filters}, i.e. those subsets $F\subseteq S$ such that, for all $a,b\in S$,
\[\tag{Filter}a,b\in F\qquad\Leftrightarrow\qquad ab\in F.\]
Note these are precisely the down-directed up-sets w.r.t. the order on $S$ defined by
\[a\leq b\qquad\Leftrightarrow\qquad a=ab\]
(as products in $S$ are meets/infima w.r.t. this ordering).  Let $\mathcal{F}(S)$ denote the \emph{filter spectrum}, i.e. space of all non-empty filters with the basis $(\mathcal{F}_a)_{a\in S}$, where
\[\mathcal{F}_a=\{F\in\mathcal{F}(S):a\in F\}.\]
Let $\mathcal{F}$ on its own denote the map $a\mapsto\mathcal{F}_a$.

\begin{prp}\label{SemilatticeFaithful}
The map $\mathcal{F}$ is a faithful spatial representation of $S$.
\end{prp}

\begin{proof}
We immediately note $\mathcal{F}_{ab}=\mathcal{F}_a\cap\mathcal{F}_b$, for all $a,b\in S$, thanks to the defining property of filters, i.e. $\mathcal{F}$ is a spatial representation.  To see that $\mathcal{F}$ is also faithful, first note that, for each $a\in S$, there is a minimal filter containing $a$, namely
\[a^\leq=\{b\in S:a\leq b\}.\]
Thus, for any $a,b\in S$, $\mathcal{F}_a=\mathcal{F}_b$ would imply $a^\leq\in\mathcal{F}_a=\mathcal{F}_b$ and $b^\leq\in\mathcal{F}_b=\mathcal{F}_a$, i.e. $b\in a^\leq$ and $a\in b^\leq$ and hence $a=ab=b$.
\end{proof}

We call $\mathcal{F}$ the \emph{filter representation} of $S$.  Note, however, that this is not the only faithful representation.  Indeed, the proof above shows that the subrepresentation obtained by restricting to principal filters (i.e. those of the form $a^\leq$) is still faithful.  What distinguishes the filter representation is its universality.

Roughly speaking, this means that the filter representation contains all possible spatial representations.  To be more precise, first note that any continuous function $\phi:Y\rightarrow X$ between topological spaces $X$ and $Y$ defines a semilattice homomorphism $\overline\phi:\mathcal{O}(X)\rightarrow\mathcal{O}(Y)$ (in the opposite direction) by taking preimages, i.e.
\[\overline\phi(O)=\phi^{-1}[O].\]
In particular, if $\theta:S\rightarrow\mathcal{O}(X)$ is a spatial representation on $X$ and $\phi:Y\rightarrow X$ is a continuous map on another space $Y$, then we immediately obtain a spatial representation $\overline\phi\circ\theta$ on $Y$ by simply composing $\theta$ with $\overline\phi$, i.e. for all $a\in S$,
\[\overline\phi\circ\theta(a)=\phi^{-1}[\theta(a)].\]
We call a spatial representation $\mu:S\rightarrow\mathcal{O}(X)$ \emph{universal} if every spatial representation $\theta:S\rightarrow\mathcal{O}(Y)$ can be written uniquely as such a composition, i.e. if there exists a unique continuous map $\phi:Y\rightarrow X$ such that $\theta=\overline\phi\circ\mu$.

\begin{prp}\label{SemilatticeUniversal}
The filter representation $\mathcal{F}$ of any semilattice $S$ is universal.
\end{prp}

\begin{proof}
Take any spatial representation $\theta:S\rightarrow\mathcal{O}(X)$.  For each $x\in X$, define
\[\phi(x)=\{a\in S:x\in\theta(a)\}.\]
Note $a,b\in\phi(x)$ iff $x\in\theta(a)\cap\theta(b)=\theta(ab)$ iff $ab\in\phi(x)$, which shows that $\phi(x)$ is a filter.  Moreover, $\phi(x)$ is also non-empty, as $\theta[S]$ covers $X$, so $\phi(x)\in\mathcal{F}(S)$.  Further note that $\phi(x)\in\mathcal{F}_a$ iff $a\in\phi(x)$ iff $x\in\theta(a)$ and hence
\[\phi^{-1}[\mathcal{F}_a]=\theta(a)\in\mathcal{O}(X).\]
Thus $\phi$ is a continuous map from $X$ to $\mathcal{F}(S)$ with $\theta=\overline\phi\circ\mathcal{F}$.

For uniqueness, say we had another map $\phi':X\rightarrow\mathcal{F}(S)$ with $\overline{\phi'}\circ\mathcal{F}=\theta$, i.e. $\phi'^{-1}[\mathcal{F}_a]=\theta(a)$, for all $a\in S$.  Then again we see that $x\in\theta(a)$ iff $\phi'(x)\in\mathcal{F}_a$ iff $a\in\phi'(x)$ and hence $\phi'(x)=\{a\in S:x\in\theta(a)\}=\phi(x)$, for all $x\in X$.
\end{proof}

\begin{rmk}
The filter representation has various other subrepresentations that are of interest, particularly if one wants to represent semilattices on nicer topological spaces, e.g. with better separation properties.

For example, one could restrict to irreducible filters, which are always plentiful enough to retain faithfulness -- see \cite[Corollary 9]{Celani2003}.  The irreducible subrepresentation can even be extended to a categorical duality (see \cite{CelaniGonzalez2020}) generalising the classic Stone duality between distributive lattices and spectral spaces (see \cite{Stone1938}).  The filter representation can also be extended to different but analogous duality between semilattices and `HMS spaces' -- see \cite{JipsenMoshier2014}.

One could also restrict to ultrafilters, i.e. maximal proper filters, as in Wallman's work (see \cite{Wallman1938}).  The space of ultrafilters is necessarily $T_1$, although it is only faithful for special kinds of semilattices, e.g. separative semilattices with $0$.  We will also briefly examine ultrafilters in more general semigroups later in \autoref{IdealUltrafilters}.  Indeed, the present paper could be viewed as a prelude to our subsequent work in \cite{Bice2021Pierce} and \cite{Bice2021DHKR} where ultrafilters play a much greater role in non-commutative extensions of the classic Gelfand duality (see \cite{Gelfand1941} and \cite{GelfandNaimark1943}).
\end{rmk}

\subsection{Inverse Semigroups}\label{InverseSemigroups}
Our goal in the present paper is to extend \autoref{SemilatticeFaithful} and \autoref{SemilatticeUniversal} to non-commutative semigroups.  As a first step in this direction, let us consider inverse semigroups (see \cite{Lawson1998}), i.e. semigroups $S$ such that each $a\in S$ has a unique generalised inverse $a^{-1}$, meaning
\[aa^{-1}a=a\qquad\text{and}\qquad a^{-1}aa^{-1}=a^{-1}.\]
So semilattices are inverse semigroups where $a^{-1}=a$.  Moreover, the idempotents
\[\mathsf{E}(S)=\{e\in S:ee=e\}\]
of any inverse semigroup $S$ always form a semilattice.  Consequently, inverse semigroups constitute a natural non-commutative generalisation of semilattices.

Topological spaces also have a natural non-commutative generalisation, namely \'etale groupoids (see \cite{Renault1980} and \cite{Sims2018}).  Let us digress for a moment to recall all the relevant definitions.  First, a \emph{groupoid} is a category where all the morphisms are isomorphisms.  Following standard practice, we generally forget about the objects and just consider groupoids as collections of morphisms on which we have an inverse operation and a partially defined product.  Given a groupoid $G$ let
\[G^0=\{g\in G:gg=g\}.\]
These are the \emph{units} of $G$.  For any $g\in G$, we denote the \emph{source} and \emph{range} units by
\[\mathsf{s}(g)=g^{-1}g\qquad\text{and}\qquad\mathsf{r}(g)=gg^{-1}.\]
These determine when the product is defined, i.e. $gh$ is defined iff $(g,h)\in G^2$ where
\[G^2=\{(g,h)\in G\times G:\mathsf{s}(g)=\mathsf{r}(h)\}.\]
As usual, we call $B\subseteq G$ a \emph{bisection} or \emph{slice} of a groupoid $G$ if the source and range maps $\mathsf{s},\mathsf{r}:G\rightarrow G^0$ are injective on $B$ or, equivalently, if $BB^{-1}\cup B^{-1}B\subseteq G^0$, where $B^{-1}=\{b^{-1}:b\in B\}$ and $AB=\{ab:a\in A\text{ and }b\in B\}$.

A \emph{topological groupoid} is a groupoid $G$ with a topology which makes the inverse $g\mapsto g^{-1}$ and product $(g,h)\mapsto gh$ continuous (on $G^2$).  A topological groupoid is \emph{\'etale} if the source $\mathsf{s}$ (and hence the range $\mathsf{r}$ and product too) is also an open map, i.e. $\mathsf{s}$ maps open subsets of $G$ to open subsets of $G$ (which, in particular, implies that the units $G^0=\mathsf{s}[G]$ form an open subset).  Equivalently, a topology makes a groupoid \'etale precisely when the open bisections
\[\mathcal{B}(G)=\{B\in\mathcal{O}(G):BB^{-1}\cup B^{-1}B\subseteq G^0\}\]
form both a basis for the topology of $G$ and an inverse semigroup under pointwise products $(A,B)\mapsto AB$ and inverses $B\mapsto B^{-1}$ (see \cite[Theorem 5.18]{Resende2007} and \cite[Proposition 6.6]{BiceStarling2018}).  This characterisation suggests that \'etale groupoids are natural structures on which to represent inverse semigroups.

Accordingly, we call a map $\theta:S\rightarrow\mathcal{B}(G)$ an \emph{\'etale representation} of an inverse semigroup $S$ if $\theta[S]$ covers the \'etale groupoid $G$ and, for all $a,b\in S$,
\[\theta(ab)=\theta(a)\theta(b)\]
(note this is consistent with our previous notion of a spatial representation, as $AB=A\cap B$ for any $A,B\subseteq G^0$, in particular for any $A,B\in\mathcal{O}(X)$ when we consider a space $X$ as an \'etale groupoid with the trivial product $xx=x$, for all $x\in X$).  To obtain a faithful \'etale representation of an inverse semigroup $S$ we can again consider filters with respect to the canonical ordering on $S$ given by
\[a\leq b\qquad\Leftrightarrow\qquad a=aa^{-1}b.\]
So, as before, we let $\mathcal{F}(S)$ denote the non-empty filters in $S$, i.e. the down-directed up-sets.  These form an \'etale groupoid with inverse $F\mapsto F^{-1}$ and product
\[E\cdot F=(EF)^\leq\quad\text{when}\quad(E^{-1}E)^\leq=(FF^{-1})^\leq,\]
where $T^\leq=\{s\geq t:t\in T\}$ and the topology again has basis $(\mathcal{F}_a)_{a\in S}$, where $\mathcal{F}_a=\{F\in\mathcal{F}(S):a\in F\}$ -- see \cite[Proposition 2.7, Example 2.12, Theorem 2.32 and Proposition 2.38]{Bice2021}, which extend results from \cite{LawsonLenz2013} and \cite{LawsonMargolisSteinberg2013}.  We call $\mathcal{F}(S)$ with this \'etale groupoid structure the \emph{filter groupoid} of $S$.  Again we let $\mathcal{F}$ on its own denote the map $a\mapsto\mathcal{F}_a$.

\begin{prp}\label{InverseSemigroupFaithful}
The map $\mathcal{F}$ is a faithful \'etale representation of $S$.
\end{prp}

\begin{proof}
For any $a,b\in S$, we immediately see that $\mathcal{F}_a\cdot\mathcal{F}_b\subseteq\mathcal{F}_{ab}$.  Conversely, for any $F\in\mathcal{F}_{ab}$, set $E=(Fb^{-1})^\leq\in\mathcal{F}_a$, as $a\geq abb^{-1}\in Fb^{-1}$.  For any $s,t\in S$ and $e\in\mathsf{E}(S)$, we see that $set\leq st$ and hence
\[(FF^{-1})^\leq\subseteq(Fb^{-1}bF^{-1})^\leq\subseteq(EE^{-1})^\leq\]
For the reverse inclusion, first note that $FF^{-1}F\subseteq F^\leq=F$ -- for any $f,g,h\in F$, we have $i\in F$ with $i\leq f,g,h$, as $F$ is a filter, and hence $i=ii^{-1}i\leq fg^{-1}h$.  Thus \[Eb=(Fb^{-1})^\leq b\subseteq(Fb^{-1}b)^\leq\subseteq F^\leq=F.\]
As above, it then follows that $(EE^{-1})^\leq\subseteq(Ebb^{-1}E^{-1})^\leq\subseteq(FF^{-1})^\leq$.  Thus $E^{-1}\cdot F$ is defined and
\[F=E\cdot E^{-1}\cdot F\in\mathcal{F}_a\cdot\mathcal{F}_{a^{-1}}\cdot\mathcal{F}_{ab}\subseteq\mathcal{F}_a\cdot\mathcal{F}_{a^{-1}ab}\subseteq\mathcal{F}_a\cdot\mathcal{F}_b.\]
Thus $\mathcal{F}_{ab}=\mathcal{F}_a\cdot\mathcal{F}_b$, showing that $a\mapsto\mathcal{F}_a$ is indeed an \'etale representation.  As before, we see the representation is faithful by considering principal filters.
\end{proof}

Of course, this is still not the only faithful \'etale representation and one would expect the filter representation to again be distinguished by its universality.  However, universality fails, at least with respect to the obvious morphisms.

Specifically, we call a map $\phi:G\rightarrow H$ between \'etale groupoids an \emph{\'etale morphism} if $\phi$ is a continuous star-bijective functor.  Here \emph{star-bijectivity} means that $\phi$ maps each star $Ge$, for $e\in G^0$, bijectively on the corresponding star $H\phi(e)$ (see \cite{Brown1970} and \cite[Ch 13]{Higgins1971} -- star-bijective functors are also called \emph{covering functors} in \cite{Lawson2012}).  As noted below in \autoref{Morphisms}, star-bijectivity ensures that preimages of slices are slices and that preimages respect products, i.e. for all $A,B\in\mathcal{B}(G)$,
\[\phi^{-1}[AB]=\phi^{-1}[A]\phi^{-1}[B].\]
Thus any \'etale morphism $\phi:H\rightarrow G$ again defines a semigroup homomorphism $\overline\phi:\mathcal{B}(G)\rightarrow\mathcal{B}(H)$ (in the opposite direction) between their open slice semigroups by just taking preimages.  Composing $\overline\phi$ with any \'etale representation $\theta:S\rightarrow\mathcal{B}(G)$ on $G$ then yields an \'etale representation $\overline\phi\circ\theta:S\rightarrow\mathcal{B}(H)$ on $H$.

Just like before, we call an \'etale representation $\mu:S\rightarrow\mathcal{B}(G)$ \emph{universal} if every \'etale representation $\theta:S\rightarrow\mathcal{B}(H)$ can be written uniquely as such a composition, i.e. there exists a unique \'etale morphism $\phi:H\rightarrow G$ such that $\theta=\overline\phi\circ\mu$.  To see why the proof of \autoref{SemilatticeUniversal} no longer works note that, while we can still define $\phi(g)=\{a\in S:g\in\theta(a)\}$, this may not be a filter anymore.  To construct universal representations here we must consider a slightly broader inverse semigroup generalisation of filters in semilattices.

Specifically, we call $C\subseteq S$ a \emph{coset} if, for all $a,b\in C$ and $c\in S$,
\[\tag{Coset}\label{Coset}c\in C\qquad\Leftrightarrow\qquad ab^{-1}c\in C.\]
Equivalently, $C\subseteq S$ is a coset when $C=C^\leq=CC^{-1}C$ (see \cite[\S1.4 Proposition 26]{Lawson1998}).  The non-empty cosets $\mathcal{C}(S)$ again form an \'etale groupoid with inverse $C\mapsto C^{-1}$, product $B\cdot C=(BC)^\leq$ and subbasis $(\mathcal{C}_a)_{a\in S}$, where
\[\mathcal{C}_a=\{C\in\mathcal{C}(S):a\in C\}\]
(again see \cite[Example 2.12 and Theorem 2.32]{Bice2021}).  We call $\mathcal{C}(S)$ with this \'etale groupoid structure the \emph{coset groupoid} of $S$ and again let $\mathcal{C}$ denote the map $a\mapsto\mathcal{C}_a$.  As in the proof of \autoref{InverseSemigroupFaithful}, we see that $\mathcal{C}$ is a faithful \'etale representation of $S$ which we call the \emph{coset representation}.

\begin{prp}
The coset representation of any inverse semigroup $S$ is universal.
\end{prp}

\begin{proof}
Take any \'etale representation $\theta:S\rightarrow\mathcal{B}(G)$.  For each $g\in G$, define
\[\phi(g)=\{a\in S:g\in\theta(a)\}.\]
If $a,b,c\in\phi(g)$ then $g=gg^{-1}g\in\theta(a)\theta(b)^{-1}\theta(c)=\theta(ab^{-1}c)$ which means that $ab^{-1}c\in\phi(g)$.  On the other hand, if $a,b,ab^{-1}c\in\phi(g)$ then
\[g=gg^{-1}g\in\theta(b)\theta(a)^{-1}\theta(ab^{-1}c)=\theta(b)\theta(a)^{-1}\theta(a)\theta(b)^{-1}\theta(c)\subseteq G^0\theta(c)\subseteq\theta(c).\]
So $\phi(g)$ is a coset, which is also non-empty, as $(\theta(a))_{a\in S}$ covers $G$, i.e. $\phi(g)\in\mathcal{C}(S)$.

Next note that $\phi(\mathsf{s}(g))=\mathsf{s}(\phi(g))$.  Indeed, we immediately see that
\[\mathsf{s}(\phi(g))=\phi(g)^{-1}\cdot\phi(g)=(\phi(g)^{-1}\phi(g))^\leq\subseteq\phi(g^{-1}g)=\phi(\mathsf{s}(g)).\]
Conversely, for any $a\in\phi(\mathsf{s}(g))$, taking $b\in\phi(g)$ we see that $a\geq ab^{-1}b\in\phi(g^{-1})\phi(g)$ and hence $a\in(\phi(g^{-1})\phi(g))^\leq=\phi(g^{-1})\cdot\phi(g)=\phi(\mathsf{s}(g))$, i.e. $\phi(\mathsf{s}(g))\subseteq\mathsf{s}(\phi(g))$.  Likewise, $\phi(\mathsf{r}(g))=\mathsf{r}(\phi(g))$ so $\mathsf{s}(g)=\mathsf{r}(h)$ implies $\mathsf{s}(\phi(g))=\mathsf{r}(\phi(h))$ and hence $\phi(g)\cdot\phi(h)$ is a valid product.  Again $\phi(g)\cdot\phi(h)\subseteq\phi(gh)$ is immediate while the reverse inclusion follows as above.  This shows that $\phi$ is functor.

To show that $\phi$ is star-bijective, take any $i\in G^0$.  If $\mathsf{s}(g)=\mathsf{s}(h)=i$ and $\phi(g)=\phi(h)$ then, for any $a\in\phi(g)=\phi(h)$, we see that $g,h\in\theta(a)$ and hence $g=h$, as $\theta(a)$ is a slice, showing that $\phi$ is star-injective.  On the other hand, for any $C\in\mathcal{C}(S)$ with $\mathsf{s}(C)=\phi(i)$, taking $c\in C$ we see that $c^{-1}c\in C^{-1}C\subseteq\mathsf{s}(C)=\phi(i)$ and hence $i\in\theta(c^{-1}c)=\theta(c)^{-1}\theta(c)=\mathsf{s}[\theta(c)]$, as $\theta(c)$ is slice.  Thus we have $g\in\theta(c)$ with $\mathsf{s}(g)=i$ and hence $\phi(g)=\phi(g)\cdot\phi(i)=(c\phi(i))^\leq=C\cdot\mathsf{s}(C)=C$, by \cite[Theorem 2.32]{Bice2021}, proving that $\phi$ is also star-surjective.

For continuity note that, as before, $\phi(g)\in\mathcal{C}_a$ iff $a\in\phi(g)$ iff $g\in\theta(a)$ and hence
\[\phi^{-1}[\mathcal{C}_a]=\theta(a)\in\mathcal{O}(G).\]
Thus $\phi$ is an \'etale morphism from $G$ to $\mathcal{C}(S)$ with $\theta=\overline\phi\circ\mathcal{C}$.  Uniqueness follows by the same argument as before.
\end{proof}

As hinted at above, there is actually an alternative approach here -- instead of broadening the class of filters, we can broaden the class of \'etale morphisms.  Specifically, we can consider more general relational Zakrzewski morphisms between \'etale groupoids -- see \autoref{Morphisms} below.  Indeed, there is a Zakrzewski morphism between filters and cosets of any inverse semigroup, which can be used to show that the filter representation is universal with respect to Zakrzewski morphisms (see \autoref{triEtaleMorphism} and \autoref{ZakrzewskiUniversalEtaleRepresentation} below for more general structured semigroup results).

\subsection{Structured Semigroups}\label{StructuredSemigroups}

The only problem with inverse semigroups is that they are not general enough, not if we want to apply the theory to semigroups commonly arising in other fields.  Our primary motivating examples here come from operator algebras, where Cartan subalgebras of C*-algebras have received increasing attention in recent years (see \cite{Renault2008}).  A key role in their theory is played by the \emph{normaliser semigroup} $S$ of a Cartan subalgebra $C\subseteq A$ given by
\[\tag{Normalisers}S=\{a\in A:aCa^*+a^*Ca\subseteq C\}.\]
Note if $A$ is a projectionless C*-algebra, like the well-known Jiang-Su algebra (which does indeed have Cartan subalgebras -- see \cite{DeeleyPutnamStrung2015}), then $S$ can not contain any non-trivial idempotents and thus certainly can not be an inverse semigroup.  Even when $A$ does have a large supply of projections, e.g. when $A$ has real rank zero, there is still no guarantee that these projections will lie in $S$.

However, the Cartan subalgebra $C$ here is a subsemigroup of $S$ which still behaves much like the idempotents in an inverse semigroup.  For one thing, $C$ is commutative, and its self-adjoint part $C_\mathrm{sa}=\{a\in C:a=a^*\}$ also has a natural lattice structure.  This exposes the possibility that inverse semigroup techniques could still be applied to Cartan pairs if we could somehow extend the relevant theory to *-semigroups $S$ with a sufficiently well-behaved *-subsemigroup $C$.  This is precisely what we did in our previous paper \cite{Bice2021}, where we examined `Weyl *-semigroups' and their coset and filter groupoids.

In the present paper we have our eye to applications even further afield (see \cite{Bice2021Pierce} and \cite{Bice2021DHKR}), namely to algebraic analogues of C*-algebras known as Steinberg algebras (see \cite{Steinberg2010} and \cite{ClarkFarthingSimsTomforde2014}) as well as to non-commutative Cartan subalgebras of C*-algebras (see \cite{Exel2011} and \cite{KwasniewskiMeyer2020}).  This forces us to consider more general semigroups $S$ without an involution $^*$, with more general non-commutative subsemigroups $C$.  However, we can always replace $C$ with its centre
\[\mathsf{Z}(C)=\{z\in C:\forall c\in C\ (cz=zc)\},\]
or some fixed subsemigroup $Z$, whenever commutativity is vital.

This naturally leads us to the concept of a `structured semigroup'.  Formally this a triple $(S,N,Z)$ where $S$ is a semigroup with distinguished subsemigroups $N$ and $Z\subseteq\mathsf{Z}(N)$ satisfying certain weak normality conditions.  To describe these conditions, first recall that $N\subseteq S$ is said to be \emph{normal} if, for all $s\in S$,
\[\tag{Normal}sN=Ns.\]
Let us call $Z\subseteq S$ \emph{binormal} if, for all $a,b\in S$,
\[\tag{Binormal}ab,ba\in Z\qquad\Rightarrow\qquad aZb\cup bZa\subseteq Z.\]
If $N$ is another subset, we call it \emph{$Z$-trinormal} if, for all $a,b\in S$ and $n\in N$,
\[\tag{$Z$-Trinormal}\label{ZTrinormal}abn=n=nab\quad\text{and}\quad ab,ba\in Z\qquad\Rightarrow\qquad bna\in N.\]
Note that if $Z\subseteq\mathsf{Z}(N)$, as in the structured semigroups considered below, then we always have $abn=nab$ whenever $ab\in Z$ and $n\in N$.  So, in this case, the left side of \eqref{ZTrinormal} could be replaced with just $abn=n$ or $nab=n$.

\begin{dfn}\label{StructuredSemigroupDef}
We call $(S,N,Z)$ a \emph{structured semigroup} if
\begin{enumerate}
\item $S$ is a semigroup.
\item $N$ is a $Z$-trinormal subsemigroup of $S$.
\item $Z\subseteq\mathsf{Z}(N)$ is a binormal subsemigroup of $S$.
\end{enumerate}
\end{dfn}

Note that every normal subsemigroup $Z$ is binormal, as $ab\in Z$ then implies that $aZb=abZ\subseteq ZZ\subseteq Z$.  However, the converse can fail, e.g. Cartan subalgebras are always binormal in their normaliser semigroups, even though they are not necessarily normal \textendash\, see \cite[Example 7.3]{BiceClark2021}.  Similarly, we also immediately see that any binormal $N$ is $N$-trinormal and hence $Z$-trinormal, for any $Z\subseteq N$.  In particular, for any commutative binormal subsemigroup $C$ of a semigroup $S$, we get a structured semigroup by simply taking $N=Z=C$.  This applies to any Cartan subalgebra $C$ in its normaliser semigroup $S$, as well as to the idempotents in any inverse semigroup, thanks to the following result.

\begin{prp}\label{NormalIdempotents}
The idempotents $\mathsf{E}(S)$ of any inverse semigroup $S$ are normal.
\end{prp}

\begin{proof}
For any $s\in S$ and $e\in\mathsf{E}(S)$, note that $ses^{-1}ses^{-1}=ss^{-1}sees^{-1}=ses^{-1}$, i.e. $ses^{-1}\mathsf{E}(S)$.  Moreover, $se=ss^{-1}se=ses^{-1}s$, showing that $s\mathsf{E}(S)\subseteq\mathsf{E}(S)s$, while the reverse inclusion follows by a dual argument.
\end{proof}

However, inverse semigroups can also form structured semigroups in other ways, i.e. there may be other natural choices for $N$ and $Z$, e.g. see \autoref{SliceSectionSemigroups} below.

The next question is how structured semigroups should be represented.  As before, we could consider \'etale representations, i.e. semigroup homomorphisms $\theta:S\rightarrow\mathcal{B}(G)$, for some \'etale groupoid $G$.  As $N$ is meant to be a kind of diagonal of $S$, it would be natural to further require that $\theta$ maps $N$ onto open subsets of the diagonal of $G$, i.e. its unit space $G^0$, as we do in \autoref{EtaleRep} below.  However, these are the idempotents of $\mathcal{B}(G)$, so if $\theta$ were also faithful then this would imply that $N$ also consists entirely of idempotents.  So again we see that these kinds of \'etale representations are not general enough, not if we want structured semigroups with few idempotents to potentially have faithful representations.

To fix this, we again take inspiration from operator algebra theory, specifically Kumjian and Renault's Weyl groupoid construction from a Cartan pair $(A,C)$ (see \cite{Kumjian1986} and \cite{Renault2008}).  This results in the C*-algebra $A$ being represented as continuous sections of a saturated Fell line bundle $\pi:F\rightarrow G$ over a locally compact \'etale groupoid $G$.  The sections coming from the normaliser semigroup $S$ are furthermore supported on open slices of $G$.  Restricting to these supports yields a representation of $S$ as partial sections defined on slices taking values in the groupoid $F^\times$ of invertible elements.  This suggests that we should likewise represent more general structured semigroups as partial sections of groupoid bundles.  Accordingly, we start the paper proper with a discussion of these bundles and their slice-sections.

\section{Slice-Sections}\label{SliceSections}

We assume the reader is familiar with the basics of \'etale groupiods, as briefly outlined in the previous section.  We refer the reader to \cite{Renault1980} and \cite{Sims2018} for further background, the key difference in our work being that our \'etale groupoids (and even their unit spaces) can be highly non-Hausdorff.

Recall that an \emph{isocofibration} $\pi:C\rightarrow D$ between categories $C$ and $D$ is a functor that is injective on objects/units or, equivalently, such that $cd$ is defined in $C$ whenever $\pi(c)\pi(d)$ is defined in $D$.  Also a map $\pi:X\rightarrow Y$ on a topological space $X$ is \emph{locally injective} if every point in $X$ has a neighbourhood on which $\pi$ is injective or, equivalently, such that the open subsets of $X$ on which $\pi$ is injective cover $X$.

\begin{dfn}
We call $\pi:F\rightarrow G$ a \emph{groupoid bundle} if
\begin{enumerate}
\item $F$ is a topological groupoid,
\item $G$ is an \'etale groupoid, and
\item $\pi$ is an open continuous isocofibration.
\end{enumerate}
An \emph{\'etale bundle} is a groupoid bundle $\pi:F\rightarrow G$ that is also locally injective.
\end{dfn}

If $\pi:F\rightarrow G$ is an \'etale bundle then, in particular, $\pi$ is an open continuous and locally injective map, i.e. a local homeomorphism.  In this case, the source map $\mathsf{s}$ on $F$ is also a local homeomorphism, which means $F$ is also an \'etale groupoid.

\begin{prp}
If $\pi:F\rightarrow G$ is an \'etale bundle then $F$ is an \'etale groupoid.
\end{prp}

\begin{proof}
Assume $\pi:F\rightarrow G$ is an \'etale bundle.  First we show that $F^0$ is open in $F$.  To see this, take any $f\in F^0$.  As $\pi$ is locally injective, $f$ has an open neighbourhood $O$ on which $\pi$ is injective.  We claim that
\[f\in O\cap\mathsf{s}^{-1}[O]\cap\pi^{-1}[G^0]\subseteq F^0.\]
To start with, $f\in F^0$ implies that $\mathsf{s}(f)=f\in O$ and $\pi(f)\in G^0$, as $\pi$ is a functor, showing that $f\in O\cap\mathsf{s}^{-1}[O]\cap\pi^{-1}[G^0]$.  On the other hand, for any $e\in O\cap\mathsf{s}^{-1}[O]\cap\pi^{-1}[G^0]$, we see that $e,\mathsf{s}(e)\in O$ and $\pi(e)=\mathsf{s}(\pi(e))=\pi(\mathsf{s}(e))$ and hence $e=\mathsf{s}(e)\in F^0$, as $\pi$ is injective on $O$, proving the claim.  As $\mathsf{s}$ and $\pi$ are continuous, and $O$ and $G^0$ are open, it follows that $O\cap\mathsf{s}^{-1}[O]\cap\pi^{-1}[G^0]$ is also open.  As $f$ was arbitrary, this shows that $F^0$ is open.

As $\pi$ is injective on units, for any $O\subseteq F$,
\[\mathsf{s}[O]=F^0\cap\pi^{-1}[\mathsf{s}[\pi[O]]].\]
If $O$ is open then this shows that $\mathsf{s}[O]$ is also open, i.e. $\mathsf{s}$ is an open map on $F$ and hence $F$ is an \'etale groupoid.
\end{proof}

The simplest examples of groupoid bundles come from the obvious structure on the Cartesian product of any \'etale groupoid and any topological group.

\begin{xpl}[Trivial Bundles]
If $G$ is an \'etale groupoid and $T$ is a topological group then $G\times T$ is a topological groupoid in the product topology where
\[(g,t)(h,u)=(gh,tu),\]
when $gh$ is defined.  The projection $\pi:G\times T\rightarrow G$ is then a groupoid bundle.  If $T$ is discrete then $\pi$ is also locally injective and hence a local homeomorphism.
\end{xpl}

More interesting `twisted' bundles can be obtained by modifying the structure above via a given $2$-cocycle, just like in \cite{Renault1980} and \cite[\S5.1]{Sims2018}.

\begin{xpl}[Twisted Bundles]
Again say $G$ is an \'etale groupoid and $T$ is a topological group.  Furthermore, say $\sigma:G^2\mapsto T$ is a continuous $2$-cocycle, i.e.
\begin{equation}\label{cocycle1}
\sigma(\mathsf{r}(g),g)=1=\sigma(g,\mathsf{s}(g)),
\end{equation}
for all $g\in G$, and, for all $t\in T$ and $g,h,i\in G$ such that $ghi$ is defined,
\begin{equation}\label{cocycle}
\sigma(g,h)t\sigma(gh,i)=\sigma(g,hi)t\sigma(h,i).
\end{equation}
We claim $G\times T$ is a topological groupoid in the product topology where this time
\[(g,t)(h,u)=(gh,t\sigma(g,h)u),\]
when $gh$ is defined.  Indeed, we verify associativity by noting that
\begin{align*}
[(g,t)(h,u)](i,v)&=(gh,t\sigma(g,h)u)(i,v)=(ghi,t\sigma(g,h)u\sigma(gh,i)v)\\
(g,t)[(h,u)(i,v)]&=(g,t)(hi,u\sigma(h,i)v)=(ghi,t\sigma(g,hi)u\sigma(h,i)v),
\end{align*}
where the last two expressions are equal by \eqref{cocycle}.  Likewise, for all $g\in G$,
\[\sigma(g,g^{-1})t=\sigma(g,g^{-1})t\sigma(\mathsf{r}(g),g)=\sigma(g,\mathsf{s}(g))t\sigma(g^{-1},g)=t\sigma(g^{-1},g).\]
It follows that
\begin{align*}
(g^{-1},\sigma(g,g^{-1})^{-1}t^{-1})(g,t)&=(\mathsf{s}(g),\sigma(g,g^{-1})^{-1}t^{-1}\sigma(g^{-1},g)t)\\
&=(\mathsf{s}(g),\sigma(g,g^{-1})^{-1}t^{-1}t\sigma(g,g^{-1}))\\
&=(\mathsf{s}(g),1)
\end{align*}
and, similarly, $(g,t)(g^{-1},\sigma(g,g^{-1})^{-1}t^{-1})=(\mathsf{r}(g),1)$, i.e.
\[(g,t)^{-1}=(g^{-1},\sigma(g,g^{-1})^{-1}t^{-1}).\]
As the cocycle is continuous, so is the product and inverse on $G\times T$.  So again $G\times T$ is a topological groupoid and $\pi:G\times T\rightarrow G$ is a groupoid bundle.
\end{xpl}

Our primary interest in groupoid bundles is that they form natural structures on which to represent semigroups as `slice-sections'.

\begin{dfn}
A \emph{slice-section} of a groupoid bundle $\pi:F\rightarrow G$ is a continuous map $a$ where $\mathrm{dom}(a)$ is an open slice of $G$ on which $\pi\circ a$ is the identity.
\end{dfn}

So if $a$ is a slice-section then $\mathrm{ran}(a)\subseteq F$ and $\pi(a(g))=g$, for all $g\in\mathrm{dom}(a)$.

\begin{prp}\label{SliceSectionSemigroups}
The family of all slice-sections $\mathcal{S}(\pi)$ of any groupoid bundle $\pi:F\rightarrow G$ is an inverse semigroup where, for all $a,b\in S$, the product is given by
\[ab(gh)=a(g)b(h)\quad\text{when }g\in\mathrm{dom}(a),\ h\in\mathrm{dom}(b)\text{ and }(g,h)\in G^2\]
$($so $\mathrm{dom}(ab)=\mathrm{dom}(a)\mathrm{dom}(b))$.  Moreover, $\mathcal{S}(\pi)$ has normal subsemigroups
\begin{align*}
\mathcal{N}(\pi)&=\{n\in\mathcal{S}(\pi):\mathrm{dom}(n)\subseteq G^0\}\quad\text{and}\\
\mathsf{E}(\mathcal{S}(\pi))&=\{z\in\mathcal{S}(\pi):\mathrm{ran}(z)\subseteq F^0\}\subseteq\mathsf{Z}(\mathcal{N}(\pi)).
\end{align*}
In particular, $(S,N,Z)=(\mathcal{S}(\pi),\mathcal{N}(\pi),\mathsf{E}(\mathcal{S}(\pi)))$ is a structured semigroup.
\end{prp}

\begin{proof}
If $gh$ is defined in $G$, for some $g,h\in G$, then $a(g)b(h)$ is also defined in $G$, as $\pi(a(g))=g$, $\pi(b(h))=h$ and $\pi$ is an isocofibration.  Also, for any $i\in G$, we have at most one pair $g\in\mathrm{dom}(a)$ and $h\in\mathrm{dom}(b)$ with $gh=i$, as these domains are slices.  Thus the above product yields a well-defined function $ab$ on the slice $\mathrm{dom}(ab)=\mathrm{dom}(a)\mathrm{dom}(b)$.  Furthermore, for any $g,h\in H$ such that $gh$ is defined,
\[\pi(ab(gh))=\pi(a(g)b(h))=\pi(a(g))\pi(b(h))=gh,\]
as $\pi$ is a functor, so again $\pi\circ ab$ is the identity on $\mathrm{dom}(ab)$.  As the product in $F$ is continuous, the function $ab$ is also continuous and hence a slice-section.  As the inverse map in $F$ is continuous, every slice-section $a$ has a unique inverse slice-section $a^{-1}(g)=a(g^{-1})^{-1}$, for all $g\in\mathrm{dom}(a)^{-1}$.  Thus the slice-sections do indeed form an inverse semigroup $S=\mathcal{S}(\pi)$.

Certainly the diagonal $N=\mathcal{N}(\pi)$ is a subsemigroup.  For any $a\in S$ and $n\in N$,
\[\mathrm{dom}(ana^{-1})\subseteq\mathrm{dom}(a)\mathrm{dom}(n)\mathrm{dom}(a)^{-1}\subseteq\mathrm{dom}(a)\mathrm{dom}(a)^{-1}\subseteq G^0,\]
as $\mathrm{dom}(n)\subseteq G^0$ and $\mathrm{dom}(a)$ is a slice.  Thus $an=ana^{-1}a\in Na$, showing that $aN\subseteq Na$, while the reverse inclusion follows by a dual argument, i.e. $aN=Na$, showing that $N$ is normal.  As the only idempotents in a groupoid are units, the idempotent sections are precisely those whose range consists of units, i.e.
\[Z=\mathsf{E}(S)=\{z\in S:\mathrm{ran}(z)\subseteq F^0\},\]
which is normal by \autoref{NormalIdempotents}.  As units commute with isotropy elements of any groupoid, it follows that $nz=zn$, for all $n\in N$ and $z\in Z$, i.e. $Z\subseteq\mathsf{Z}(N)$. 
\end{proof}

Note here that $N$ will be strictly bigger than $Z=\mathsf{E}(S)$, as long as there are at least some slice-sections of $\pi$ taking non-unit values.  So in this case, the above structured semigroup is different from the canonical one mentioned after \autoref{StructuredSemigroupDef} formed by simply taking $N=Z=\mathsf{E}(S)$.

We can also consider certain inverse subsemigroups of slice-sections.

\begin{xpl}
If $G$ is an \'etale groupoid with Hausdorff unit space $G^0$ then pointwise products $KL$ of compact subsets $K,L\subseteq G$ are again compact.  If $\pi:F\rightarrow G$ is a groupoid bundle then it follows that the slice-sections defined on (open) compact slices form an inverse subsemigroup $S\subseteq\mathcal{S}(\pi)$.  Here it is natural to restrict to ample groupoids (where $G^0$ and hence $G$ has a basis of compact open subsets) so these compact-slice-sections can distinguish points of $G$.

If we have a trivial ample groupoid bundle $\pi:G\times(\mathbb{F}\setminus\{0\})\rightarrow G$, where $\mathbb{F}$ is a discrete field, then any compact-slice-section $s$ can be extended to the entirety of $G$ with $0$ values outside the original domain.  The linear span of these then forms an algebra $A$ of functions from $G$ to $\mathbb{F}$ commonly known as the `Steinberg algebra' over $G$.  The compact-slice-sections of non-trivial locally injective ample groupoid bundles $\pi:F\rightarrow G$ could thus be considered as `twisted Steinberg semigroups'.
\end{xpl}

Concrete (sub)semigroups of slice-sections arising in other contexts often only have `local inverses'.

\begin{dfn}\label{local}
Given a groupoid bundle $\pi:F\rightarrow G$ and semigroup $S\subseteq\mathcal{S}(\pi)$, let
\[S_g=\{a\in S:g\in\mathrm{dom}(a)\}.\]
We call $S$ a \emph{local-inverse semigroup} if, for all $g\in G$ and $a\in S_g$, we have $a'\in S_{g^{-1}}$ and $b\in S_g$ with $a'(h^{-1})=a(h)^{-1}$, for all $h\in\mathrm{dom}(b)$.
\end{dfn}

So the values of $a'$ above are the inverses of the values of $a$ on some neighbourhood of $g$ which is also the domain of some $b\in S$.  In particular, $\mathrm{dom}(b)\subseteq\mathrm{dom}(a)$.

\begin{xpl}\label{FellBundles}
As in \cite{Kumjian1998}, we can consider a Fell bundle $\pi:F\rightarrow G$ over a locally compact \'etale groupoid $G$ with Hausdorff unit space $G^0$.  This is not quite a groupoid bundle \textendash\, even with line bundles, where the fibres are $1$-dimensional, each fibre still has a $0$.  However we can simply restrict $\pi$ to the invertible elements $F^\times$.  The slice-sections $a$ defined on open slices such that $\|a\|$ vanishes at infinity again form a subsemigroup $S\subseteq\mathcal{S}(\pi|_{F^\times})$.  Again we can extend any such slice-section to the entirety of $G$ with $0$ values outside the original domain.  The linear span of these again forms an algebra, this time with a natural reduced algebra norm.  Taking the completion then yields a C*-algebra $A$.  Thus $S$ again forms a natural semigroup from which to generate $A$.

Note the inverse of a function vanishing at infinity may not vanish at infinity and could even be unbounded, i.e. $a\in S$ does not imply $a^{-1}\in S$ so $S$ is not an inverse semigroup.  However, if $a\in S_g$ then we will have an open neighbourhood $O\subseteq\mathrm{dom}(a)$ of $g$ on which $a^{-1}$ is bounded.  We can then define a function $a'\in S_{g^{-1}}$ which vanishes at infinity but coincides with $a^{-1}$ on $O^{-1}$.  Taking any $b\in S_g$ with $\mathrm{dom}(b)\subseteq O$ then witnesses the fact that $S$ is a local-inverse semigroup.

Similarly, the sections of $\pi$ with range in $F^0$ will rarely vanish at infinity (only if their domain is also compact).  Thus it is more natural to allow scalar multiples of units as well, i.e. instead of taking $Z$ to be the idempotent slice-sections as in \autoref{SliceSectionSemigroups}, we can take $Z=S\cap\mathcal{Z}(\pi)$, where $\mathcal{Z}(\pi)$ denotes the \emph{central-diagonal}
\begin{equation}\label{ScalarZ}
\mathcal{Z}(\pi)=\{z\in\mathcal{S}(\pi):\mathrm{ran}(z)\subseteq\mathbb{C}F^0\}.
\end{equation}
As the name suggests, $\mathcal{Z}(\pi)$ is a central subsemigroup of the diagonal $\mathcal{N}(\pi)$ and hence $Z=S\cap\mathcal{Z}(\pi)$ is also a central subsemigroup of $N=S\cap\mathcal{N}(\pi)$.  Moreover, note that we can always choose the $a'$ in the previous paragraph so that $aa',a'a\in Z$.  Also, unlike in the previous examples, $N$ and $Z$ here may not be normal (e.g. in \cite[Example 7.3]{BiceClark2021} mentioned above) but they will still be binormal and hence $(S,N,Z)$ will still form a structured semigroup.

Another possibility is to leave $\pi$ as is \textendash\, even though $F$ is only a category, not a groupoid, the slice-sections $\mathcal{S}(\pi)$ still form a semigroup with a subsemigroup $S$ of functions vanishing at infinity.  Again taking $N=S\cap\mathcal{N}(\pi)$ and $Z=S\cap\mathcal{Z}(\pi)$, we see that $N$ may not even be binormal now, however it will still be $Z$-trinormal, i.e. $(S,N,Z)$ will still form a structured semigroup.
\end{xpl}

Our main goal will be to obtain a kind of converse to \autoref{SliceSectionSemigroups}.  More precisely, we want to show how to represent a structured semigroup as a subsemigroup of $\mathcal{S}(\pi)$, for some \'etale bundle $\pi:F\rightarrow G$, which is universal for an appropriate class of representations on even more general groupoid bundles.

To get some hint as to how we might do this, take an \'etale bundle $\pi:F\rightarrow G$ and consider the inverse semigroup of all slice-sections $S=\mathcal{S}(\pi)$.  The sets of the form $S_g$, for $g\in G$, form `cosets' in an appropriate sense, generalising the usual cosets one considers in inverse semigroups -- see \eqref{Coset} above.  Indeed, these general cosets can still be defined via a \emph{domination} relation $<$ on $S$ generalising the usual order on an inverse semigroup (basically $<$ corresponds to inclusion of domains \textendash\, see \autoref{<Equiv} below).  Moreover, for any $g,h\in G$ such that $gh$ is defined, the products $S_gS_h$ of slice-sections in $S_g$ and $S_h$ are $<$-coinitial in the family of slice-sections $S_{gh}$.  In this way, the product structure of $S$ encodes that of $G$.

For each $f\in F$ and $g\in G$ with $\pi(f)=g$, we can further consider the subfamily $S^f\subseteq S_g$ of slice-sections taking the value $f$ at $g$
\[S^f=\{a\in S:a(g)=f\}.\]
These form equivalence classes in $S_g$ modulo a relation $\sim_{S_g}$ again determined by the structured semigroup.  Moreover, for any $e,f\in F$ such that $ef$ is defined, $S^{ef}$ coincides with the equivalence class containing $S^eS^f$, again showing how the product structure of $F$ is encoded by that of $S$.

The goal is thus to examine the domination relation $<$ in structured semigroups and the cosets and equivalence classes it gives rise to, in order to construct \'etale bundles on which to obtain universal slice-section representations.  For universality, however, we will first need to consider morphisms between groupoid bundles which correspond to semigroup homomorphisms of the associated slice-sections.

\section{Morphisms}\label{Morphisms}

Before considering groupoids bundles, we first we need an appropriate notion of morphism for \'etale groupoids.  One option would be to consider the (functional) \'etale morphisms from \autoref{InverseSemigroups}, which would be fine for \autoref{TheUniversalRepresentation} where we consider general cosets.  However, when we restrict our attention to directed cosets in \autoref{Filters}, we will need more general (relational) Zakrzewski morphisms -- see \autoref{triEtaleMorphism} below.

First let us consider any $\phi\subseteq G\times H$ as a relation `from $H$ to $G$', where $g\mathrel{\phi}h$ means $(g,h)\in\phi$.  The flip of $\phi$ is denoted by $\phi^{-1}\subseteq H\times G$, i.e. $h\mathrel{\phi^{-1}}g$ iff $g\mathrel{\phi}h$.  The \emph{image} of $H'\subseteq H$ and \emph{preimage} of $G'\subseteq G$ are denoted by
\begin{align}
\tag{Image}\phi[H']&=\{g\in G:\exists h\in H'\ (g\mathrel{\phi}h)\}.\\
\tag{Preimage}\phi^{-1}[G']&=\{h\in H:\exists g\in G'\ (g\mathrel{\phi}h)\}.
\end{align}
In particular, we define the domain and range of $\phi\subseteq G\times H$ by
\[\mathrm{dom}(\phi)=\phi^{-1}[G]\subseteq H\qquad\text{and}\qquad\mathrm{ran}(\phi)=\phi[H]\subseteq G.\]
The composition of $\phi\subseteq G\times H$ and $\psi\subseteq H\times I$ is the relation $\phi\circ\psi\subseteq G\times I$ where
\[g\mathrel{(\phi\circ\psi)}i\qquad\Leftrightarrow\qquad\exists h\in H\ (g\mathrel{\phi}h\mathrel{\psi}i).\]

\begin{dfn}
Let $G$ and $H$ be groupoids.  We call $\phi\subseteq G\times H$ \emph{functorial} if
\begin{align*}
g\mathrel{\phi}h\quad&\Rightarrow\quad g^{-1}\mathrel{\phi}h^{-1},\quad\text{and}\\
g\mathrel{\phi}h,\ g'\mathrel{\phi}h'\text{ and }\mathsf{s}(h)=\mathsf{r}(h')\quad&\Rightarrow\quad\mathsf{s}(g)=\mathsf{r}(g')\text{ and }gg'\mathrel{\phi}hh'.
\end{align*}
\end{dfn}

In particular, note $\mathrm{dom}(\phi)$ is a subgroupoid when $\phi$ is functorial.

As usual, we call $\phi\subseteq G\times H$ a \emph{function} when the image $\phi\{h\}=\phi[\{h\}]$ of every singleton $h$ in $\mathrm{dom}(\phi)$ is again a singleton in $\mathrm{ran}(\phi)$, in which case
\[g\mathrel{\phi}h\qquad\Leftrightarrow\qquad g=\phi(h).\]
So a functorial function $\phi\subseteq G\times H$ is just a functor on $\mathrm{dom}(\phi)$ in the usual sense.  We also get functions when we restrict functorial relations to slices in the range.

\begin{prp}
If $\phi\subseteq G\times H$ is functorial then
\begin{equation}\label{SliceRestriction}
B\subseteq G\text{ is a slice}\qquad\Rightarrow\qquad\phi\cap(B\times H)\text{ is function}.
\end{equation}
\end{prp}

\begin{proof}
It suffices to show that
\[f,g\mathrel{\phi}h\qquad\Rightarrow\qquad\mathsf{s}(f)=\mathsf{s}(g).\]
To see this, note that $f,g\mathrel{\phi}h$ implies $g^{-1}\mathrel{\phi}h^{-1}$, by functoriality.  As $\mathsf{s}(h)=\mathsf{r}(h^{-1})$, functoriality again yields $\mathsf{s}(f)=\mathsf{r}(g^{-1})=\mathsf{s}(g)$.
\end{proof}

We extend the terminology in \cite{Brown1970} and \cite{Higgins1971} from functions to relations.

\begin{dfn}
Assume $G$ and $H$ are groupoids and $\phi\subseteq G\times H$ is functorial.  We call $\phi$ \emph{star-injective/surjective/bijective} if, whenever $h\in H$ and $\mathsf{r}(g)\mathrel{\phi}h\in H^0$, there is at most/at least/precisely one $i\in H$ such that $\mathsf{r}(i)=h$ and $g\mathrel{\phi}i$.
\end{dfn}

More symbolically, these definitions can be summarised as follows.
\begin{align*}
\tag{Star-Injective}\mathsf{r}(g)\mathrel{\phi}h\in H^0\qquad&\Rightarrow\qquad|\phi^{-1}\{g\}\cap\mathsf{r}^{-1}\{h\}|\leq1.\\
\tag{Star-Surjective}\mathsf{r}(g)\mathrel{\phi}h\in H^0\qquad&\Rightarrow\qquad|\phi^{-1}\{g\}\cap\mathsf{r}^{-1}\{h\}|\geq1.\\
\tag{Star-Bijective}\mathsf{r}(g)\mathrel{\phi}h\in H^0\qquad&\Rightarrow\qquad|\phi^{-1}\{g\}\cap\mathsf{r}^{-1}\{h\}|=1.
\end{align*}
These have several equivalent characterisations, e.g. $\phi$ is star-injective iff its kernel is discrete, i.e. unit-preimages are units.  More important for us is the fact that star-injectivity means $\phi^{-1}$ respects slices.

\begin{prp}
Functorial $\phi\subseteq G\times H$ is star-injective iff $\phi^{-1}[G^0]\subseteq H^0$ iff
\begin{equation}\label{StarInjective}
B\text{ is a slice of }G\qquad\Rightarrow\qquad\phi^{-1}[B]\text{ is a slice of }H,
\end{equation}
\end{prp}

\begin{proof}
If $\phi$ is not star-injective then we have $h\in H^0$ and distinct $i,j\in H$ with $\mathsf{r}(h)=\mathsf{r}(i)=h$ and $g\mathrel{\phi}i,j$.  Then $G^0\ni\mathsf{s}(g)=g^{-1}g\mathrel{\phi}i^{-1}j\notin H^0$, by functoriality.

Conversely, if we have $G^0\ni g\mathrel{\phi}h\notin H^0$ then $g=\mathsf{r}(g)\mathrel{\phi}\mathsf{r}(h)\neq h$, even though $\mathsf{r}(\mathsf{r}(h))=\mathsf{r}(h)$, showing that $\phi$ is not star-injective.

Now if $\phi^{-1}[G^0]\subseteq H^0$ then, for any slice $B\subseteq G^0$, functoriality yields
\[\phi^{-1}[B](\phi^{-1}[B])^{-1}=\phi^{-1}[B]\phi^{-1}[B^{-1}]\subseteq\phi^{-1}[BB^{-1}]\subseteq\phi^{-1}[G^0]\subseteq H^0.\]
Likewise, $(\phi^{-1}[B])^{-1}\phi^{-1}[B]\subseteq H^0$, showing that $\phi^{-1}[B]$ is a slice.

Conversely, if $\phi^{-1}\{g\}$ is a slice then, in particular, $i,j\in\phi^{-1}\{g\}\cap\mathsf{r}^{-1}\{h\}$ implies $i=j$.  So if this holds for all $g\in G$ then $\phi$ is star-injective.
\end{proof}

On the other hand, star-surjectivity means $\phi^{-1}$ respects products.

\begin{prp}\label{StarSurjectiveChar}
Functorial $\phi\subseteq G\times H$ is star-surjective iff, for all $A,B\subseteq G$,
\begin{equation}\label{StarSurjective}
\phi^{-1}[AB]=\phi^{-1}[A]\phi^{-1}[B].
\end{equation}
\end{prp}

\begin{proof}
Assume $\phi\subseteq G\times H$ is functorial and star-surjective.  If $A,B\subseteq G$ then $\phi^{-1}[A]\phi^{-1}[B]\subseteq\phi^{-1}[AB]$, by the functoriality of $\phi$.  Conversely, take $a\in A$, $b\in B$ and $h\in H$ with $ab\mathrel{\phi}h$.  By functoriality, $\mathsf{r}(a)=\mathsf{r}(ab)\mathrel{\phi}\mathsf{r}(h)$.  Star-surjectivity then yields $i\in\phi^{-1}\{a\}$ with $\mathsf{r}(h)=\mathsf{r}(i)$.  Again funtoriality yields $b=a^{-1}(ab)\mathrel{\phi}i^{-1}h$ and hence $h=i(i^{-1}h)\in\phi^{-1}[A]\phi^{-1}[B]$, showing $\phi^{-1}[AB]\subseteq\phi^{-1}[A]\phi^{-1}[B]$ too.

On the other hand, if \eqref{StarSurjective} holds and $\mathsf{r}(g)\mathrel{\phi}h\in H^0$ then
\[h\in\phi^{-1}\{\mathsf{r}(g)\}=\phi^{-1}\{gg^{-1}\}=\phi^{-1}\{g\}\phi^{-1}\{g^{-1}\}=\phi^{-1}\{g\}(\phi^{-1}\{g\})^{-1}.\]
Thus we have $i\in\phi^{-1}\{g\}$ and $j\in\phi^{-1}\{g^{-1}\}$ with $h=ij$.  As $h\in H^0$, this means that $h=\mathsf{r}(h)=\mathsf{r}(ij)=\mathsf{r}(i)$, showing that $\phi$ is star-surjective.
\end{proof}

In fact, the proofs show it suffices to consider singleton $A$ and $B$ in \eqref{StarInjective} and \eqref{StarSurjective}.  We also note $\phi$ star-surjective iff $\mathrm{ran}(\phi)$ is an ideal of $G$ \textendash\, see \eqref{Ideal} below.

Extending the usual notion for functions, when $G$ and $H$ are topological spaces and $\phi\subseteq G\times H$, we call $\phi$ \emph{continuous} if $\phi^{-1}[O]$ is open, for all open $O\subseteq G$.

\begin{dfn}
A \emph{Zakrzewski morphism} $\phi\subseteq G\times H$ bewteen \'etale groupoids $G$ and $H$ is a continuous star-bijective functorial relation.
\end{dfn}

The idea of considering relational morphisms between groupoids comes from \cite{Zakrzewski1990} and was further studied in \cite{Stachura2018}.

As in \autoref{InverseSemigroups}, a functional Zakrzewski morphism will be called an \emph{\'etale morphism}.

\begin{rmk}\label{SpacesVsGroups}
Zakrzewski morphisms include both continuous functions and group homomorphisms but in opposite directions.  Indeed, if $G$ and $H$ are \'etale groupoids with $G=G^0$ and $H=H^0$, i.e. if $G$ and $H$ are just topological spaces under the trivial product defined on the diagonal, a Zakrzewski morphism $\phi\subseteq G\times H$ is just a continuous function from an open subset of $H$ to $G$.  At the other extreme, if $G^0$ and $H^0$ are singletons, i.e. if $G$ and $H$ are discrete groups, then non-empty $\phi\subseteq G\times H$ is a Zakrzewski morphism iff $\phi^{-1}\subseteq H\times G$ is a group homomorphism from $G$ to $H$.

Zakrzewski morphisms can also be seen as generalisations of groupoid bundles.  Indeed, if $F$ and $G$ are \'etale groupoids then a function $\pi:F\rightarrow G$ is a groupoid bundle iff $\pi^{-1}\subseteq F\times G$ is an open Zakrzewski morphism, i.e. a Zakrzewski morphism such that $\pi^{-1}[O]$ is open, for all open $O\subseteq G$.
\end{rmk}

Now say we have groupoid bundle $\pi:F\rightarrow G$ and a Zakrzewski morphism $\phi\subseteq G\times H$.  Consider the subspace of $F\times H$ given by
\[\phi^\pi F=\{(f,h)\in F\times H:\pi(f)\mathrel{\phi}h\}.\]
Note $\phi^\pi F$ is a topological groupoid under the product
\[(f,h)(f',h')=(ff',hh')\quad\text{when $hh'$ is defined}.\]
Let $\pi_\phi:\phi^\pi F\rightarrow H$ denote the projection onto $H$, i.e. $\pi_\phi(f,h)=h$.

\begin{prp}\label{PullbackEtale}
If $\pi:F\rightarrow G$ is a groupoid bundle and $\phi\subseteq G\times H$ is a Zakrzewski morphism then the projection $\pi_\phi:\phi^\pi F\rightarrow H$ is also a groupoid bundle.
\end{prp}

\begin{proof}
The definition of the topology and product on $\phi^\pi F$ ensures that $\pi_\phi$ is a continuous isocofibration.  It only remains to show that $\pi_\phi$ is also an open map.  To see this, take open $O\subseteq F$ and $N\subseteq H$ and note that
\begin{align*}
\pi_\phi[(O\times N)\cap\phi^\pi F]&=\{h\in N:\exists f\in O\ (\pi(f)\mathrel{\phi}h)\}\\
&=N\cap\phi^{-1}[\pi[O]],
\end{align*}
which is open because $\phi$ is continuous and $\pi$ is an open map.
\end{proof}

This $\pi_\phi:\phi^\pi F\rightarrow H$ is the \emph{pullback bundle} of $\pi$ induced by $\phi$.

\begin{dfn}
Let $\pi:F\rightarrow G$ and $\pi':F'\rightarrow G'$ be groupoid bundles.  If
\begin{enumerate}
\item $\phi\subseteq G\times G'$ is a Zakrzewski morphism,
\item $\tau:\phi^\pi F\rightarrow F'$ is a continuous functor and
\item $\pi_\phi=\pi'\circ\tau$ (i.e. $\pi'(\tau(f,g'))=g'$, for all $(f,g')\in\phi^\pi F$)
\end{enumerate}
then we call the pair $(\phi,\tau)$ a \emph{Zakrzewski-Pierce morphism} from $\pi$ to $\pi'$.\\
If $\phi$ is also a function then $(\phi,\tau)$ is a \emph{Pierce morphism}.
\end{dfn}

Similar morphisms for ring bundles were considered in \cite[Definition 6.1]{Pierce1967}, hence the name (for analogous C*-bundle morphisms, see \cite[Definition 4.3]{Varela1974}).  Like in \cite[Lemma 6.2]{Pierce1967}, Zakrzewski-Pierce morphisms of groupoid bundles naturally yield semigroup homomorphisms of their slice-sections.

\begin{thm}\label{BundleMorphism->SemigroupHomo}
If $(\phi,\tau)$ is a Zakrzewski-Pierce morphism from $\pi:F\rightarrow G$ to $\pi':F'\rightarrow G'$, we have a semigroup homomorphism $\frac{\tau}{\phi}:\mathcal{S}(\pi)\rightarrow\mathcal{S}(\pi')$ given by
\[\tfrac{\tau}{\phi}(a)=\tau\circ(((a\circ\phi)\times\mathrm{id})\circ\delta),\]
where $\delta(g')=(g',g')$ and $\mathrm{id}(g')=g'$.
\end{thm}

\begin{proof}
Assume $a$ is a slice-section of $\pi$.  In particular, $\mathrm{dom}(a)$ is an open slice, as is
\[\mathrm{dom}(\tfrac{\tau}{\phi}(a))=\mathrm{dom}(a\circ\phi)=\phi^{-1}[\mathrm{dom}(a)],\]
because $\phi$ is continuous and star-injective \textendash\, see \eqref{StarInjective}.  As $a$ and $\tau$ are also continuous, so is $\tfrac{\tau}{\phi}(a)$.  Moreover, $\phi\cap(\mathrm{dom}(a)\times G')$ and hence $a\circ\phi$ is function, by \eqref{SliceRestriction}.  Thus $\tfrac{\tau}{\phi}(a)$ is a function too such that, for all $g'\in\mathrm{dom}(\tfrac{\tau}{\phi}(a))$,
\[\tfrac{\tau}{\phi}(a)(g')=\tau((a\circ\phi)(g'),g')\]
(which is defined because $((a\circ\phi)(g'),g')\in\phi^\pi$, as $a$ is a section of $\pi$).  It follows that $\pi'(\tfrac{\tau}{\phi}(a)(g'))=\pi'(\tau((a\circ\phi)(g'),g'))=\pi_\phi((a\circ\phi)(g'),g')=g'$.  This all shows that $\tfrac{\tau}{\phi}(a)$ is a slice-section of $\pi'$, i.e. $\tfrac{\tau}{\phi}(a)\in\mathcal{S}(\pi')$.

Now given another slice-section $b\in\mathcal{S}(\pi)$, note that
\begin{align*}
\mathrm{dom}(\tfrac{\tau}{\phi}(ab))=\phi^{-1}[\mathrm{dom}(ab)]&=\phi^{-1}[\mathrm{dom}(a)\mathrm{dom}(b)],\text{ and}\\
\mathrm{dom}(\tfrac{\tau}{\phi}(a)\tfrac{\tau}{\phi}(b))=\mathrm{dom}(\tfrac{\tau}{\phi}(a))\mathrm{dom}(\tfrac{\tau}{\phi}(b))&=\phi^{-1}[\mathrm{dom}(a)]\phi^{-1}[\mathrm{dom}(b)],
\end{align*}
so $\mathrm{dom}(\tfrac{\tau}{\phi}(ab))=\mathrm{dom}(\tfrac{\tau}{\phi}(a)\tfrac{\tau}{\phi}(b))$, by the star-surjectivity of $\phi$.  Moreover, for any $g'\in\mathrm{dom}(\tfrac{\tau}{\phi}(a))$ and $h'\in\mathrm{dom}(\tfrac{\tau}{\phi}(b))$ such that $g'h'$ is defined, note that
\begin{align*}
\tfrac{\tau}{\phi}(ab)(g'h')&=\tau((ab\circ\phi)(g'h'),g'h')\\
&=\tau((a\circ\phi)(g')(b\circ\phi)(h'),g'h')\\
&=\tau(((a\circ\phi)(g'),g')((b\circ\phi)(h'),h'))\\
&=\tau((a\circ\phi)(g'),g')\tau((b\circ\phi)(h'),h')\\
&=\tfrac{\tau}{\phi}(a)(g')\tfrac{\tau}{\phi}(b)(h')\\
&=(\tfrac{\tau}{\phi}(a)\tfrac{\tau}{\phi}(b))(g'h').
\end{align*}
This shows that $\tfrac{\tau}{\phi}(ab)=\tfrac{\tau}{\phi}(a)\tfrac{\tau}{\phi}(b)$, i.e. $\tfrac{\tau}{\phi}$ is semigroup homomorphism.
\end{proof}

\section{Domination}\label{Domination}

We now make the following standing assumption throughout.

\begin{center}
$(S,N,Z)$\textbf{ is a structured semigroup (see \autoref{StructuredSemigroupDef}).}
\end{center}

\begin{dfn}
We define relations on $S$ as follows.
\begin{align*}
a<_sb\qquad&\Leftrightarrow\qquad asb=a=bsa,\ as,sa\in N\text{ and }bs,sb\in Z.\\
a<b\qquad&\Leftrightarrow\qquad\exists s\in S\ (a<_sb).
\end{align*}
When $a<_sb$, we say that \emph{$b$ dominates $a$ via $s$}
\end{dfn}

To get some intuition for domination, we imagine that $S$ is a subsemigroup of the slice-sections $\mathcal{S}(\pi)$ of some groupoid bundle $\pi:F\rightarrow G$.  As long as $N$ lies in the diagonal $\mathcal{N}(\pi)$, $a<_sb$ means that $s$ is an inverse of $b$ on the domain of $a$.

\begin{prp}\label{<Equiv}
If $\pi:F\rightarrow G$ is a groupoid bundle then, for any $a,s,b\in\mathcal{S}(\pi)$,
\[a=asb\quad\text{and}\quad as\in\mathcal{N}(\pi)\qquad\Leftrightarrow\qquad\forall g\in\mathrm{dom}(a)\ (s(g^{-1})=b(g)^{-1})\]
$($implicit in the statement $s(g^{-1})=b(g)^{-1}$ is that $g^{-1}\in\mathrm{dom}(s)$ and $g\in\mathrm{dom}(b))$.
\end{prp}

\begin{proof}
If $a=asb$ and $g\in\mathrm{dom}(a)$ then $a(g)=asb(g)=a(h)s(i)b(j)$, for some $h,i,j\in G$ with $g=hij$.  Thus $\mathsf{r}(g)=\mathsf{r}(hij)=\mathsf{r}(h)$ so $g=h$, as $\mathrm{dom}(a)$ is a slice.  If $as\in\mathcal{N}(\pi)$ too then $gi=hi\in\mathrm{dom}(as)\subseteq G^0$ and hence $i=g^{-1}$.  Thus $g=hij=gg^{-1}j=j$ so $a(g)=a(g)s(g^{-1})b(g)$ and hence $s(g^{-1})b(g)\in F^0$, i.e. $s(g^{-1})=b(g)^{-1}$.  This proves $\Rightarrow$.

Now say $s(g^{-1})=b(g)^{-1}$, for all $g\in\mathrm{dom}(a)$.  In particular, $\mathrm{dom}(a)\subseteq\mathrm{dom}(s)^{-1}$ and hence $\mathrm{dom}(as)=\mathrm{dom}(a)\mathrm{dom}(s)\subseteq\mathrm{dom}(s)^{-1}\mathrm{dom}(s)\subseteq G^0$ so $as\in\mathcal{N}(\pi)$.  Also
\[asb(g)=a(g)s(g^{-1})b(g)=a(g)b(g)^{-1}b(g)=a(g),\]
for all $g\in\mathrm{dom}(a)$.  In particular, $\mathsf{s}[\mathrm{dom}(a)]\subseteq\mathrm{dom}(sb)$ and hence
\[\mathrm{dom}(asb)=\mathrm{dom}(a)\mathrm{dom}(sb)=\mathrm{dom}(a)\]
which, together with the above computation, yields $a=asb$, proving $\Leftarrow$.
\end{proof}

In particular, if $\pi:F\rightarrow G$ is a groupoid bundle and $(S,N,Z)$ is a structured semigroup with $S\subseteq\mathcal{S}(\pi)$ and $N\subseteq\mathcal{N}(\pi)$ then it follows from the above that
\[a<b\qquad\Rightarrow\qquad\mathrm{dom}(a)\subseteq\mathrm{dom}(b).\]
On the other hand, if $S=\mathcal{S}(\pi)$, $\mathcal{N}(\pi)\subseteq N$ and $\mathsf{E}(S)\subseteq Z$ then the converse holds.  If we instead consider functions vanishing at infinity like in \autoref{FellBundles} then $<$ is instead equivalent to compact containment $\Subset$ of domains, i.e. $a<b$ means $\mathrm{dom}(a)\subseteq K\subseteq\mathrm{dom}(b)$, for some compact $K\subseteq G$ (see \cite[Proposition 4.4]{BiceClark2021}).

On inverse semigroups, $<$ is just the usual ordering.

\begin{prp}
If $S$ is an inverse semigroup and $N=Z=\mathsf{E}(S)$ then
\begin{equation}\label{CanonicalOrder}
a<b\qquad\Leftrightarrow\qquad a\in Nb\qquad\Leftrightarrow\qquad a<_{b^{-1}}b.
\end{equation}
\end{prp}

\begin{proof} If $a<b$ then we have $s\in S$ with $a=asb\in Nb$.  Conversely, if $a=nb$, for some $n\in N$, then $ab^{-1}=nbb^{-1}\in N$ and $ab^{-1}b=nbb^{-1}b=nb=a$.  Likewise $b^{-1}a=b^{-1}nb\in N$ and $bb^{-1}a=bb^{-1}nb=nbb^{-1}b=nb=a$.  As $bb^{-1},b^{-1}b\in N$, this shows that $a<_{b^{-1}}b$.
\end{proof}

Our primary goal in this section is to show that the domination relation $<$ on general structured semigroups still has many of the same properties as the usual order relation on inverse semigroups.

First we note that $<$ could have been defined in a more one-sided way.  Note here and elsewhere we will make use of the fact that a result proved from $S$ immediately yields a dual result for the opposite semigroup $S^\mathrm{op}$ (where $a\cdot^{\mathrm{op}}b=ba$).

\begin{prp}\label{1sided<}
For any $a,b,b'\in S$,
\begin{align}
a<_sb\qquad&\Leftrightarrow\qquad asb=a,\ as\in N\text{ and }bs,sb\in Z.\\
&\Leftrightarrow\qquad bsa=a,\ sa\in N\text{ and }bs,sb\in Z.
\end{align}
\end{prp}

\begin{proof}
If $asb=a$, $as\in N$ and $bs,sb\in Z$ then $as\in N$ commutes with $bs\in Z$.  Thus $bsa=bsasb=asbsb=a$ and $bsas=asbs=as\in N$ so $sa=sasb\in N$, as $N$ is $Z$-trinormal.  Thus $a<_sb$, proving the first $\Leftrightarrow$, while the second follows dually.
\end{proof}

From now on we will use the above characterisations of $<$.

\begin{prp}
For all $a,b',b,c',c\in S$,
\[\tag{Transitivity}\label{Transitivity}a<_{b'}b<_{c'}c\qquad\Rightarrow\qquad a<_{c'}c.\]
\end{prp}

\begin{proof}
If $a<_{b'}b<_{c'}c$, $ac'=ab'bc'\in NN\subseteq N$ and $ac'c=ab'bc'c=ab'b=a$.
\end{proof}

We can also switch the subscript with the right argument as follows.

\begin{prp}
For all $a,b,c,c'\in S$,
\[\tag{Switch}\label{Switch}ab\in N\quad\text{and}\quad b<_{c'}c\qquad\Rightarrow\qquad abc'<_{c}c'.\]
\end{prp}

\begin{proof}
If $ab\in N$ and $b<_{c'}c$ then $abc'c\in NZ\subseteq N$ and $abc'cc'=abc'$.
\end{proof}

Next we show that $<$ preserves the product.

\begin{prp}
For any $a,b',b,c,d',d\in S$,
\[\tag{Multiplicativity}\label{Multiplicativity}a<_{b'}b\quad\text{and}\quad c<_{d'}d\qquad\Rightarrow\qquad ac<_{d'b'}bd.\]
\end{prp}

\begin{proof}
If $a<_{b'}b$ and $c<_{d'}d$ then $b'bb'acd'=b'acd'\in NN\subseteq N$ and hence $acd'b'=bb'acd'b'\in N$, as $N$ is $Z$-trinormal.  Also $acd'b'bd=ab'bcd'd=ac$, as $b'b\in Z$ commutes with $cd'\in N$.  Moreover, as $Z$ is binormal,
\[d'b'bd\in d'Zd\subseteq Z\supseteq bZb'\ni bdd'b'.\qedhere\]
\end{proof}

Next we note multiplying by elements of $N$ on the left does not affect $<$.

\begin{prp}
For any $a,b',b\in S$ and $n\in N$,
\[\tag{$N$-Invariance}\label{NInvariance}a<_{b'}b\qquad\Rightarrow\qquad an,na<_{b'}b.\]
\end{prp}

\begin{proof}
If $a<_{b'}b$ then $nab'\in NN\subseteq N$ and $nab'b=na$ and hence $na<_{b'}b$, while $an<_{b'}b$ follows by duality.
\end{proof}

We also have a similar result for elements of $Z$.

\begin{prp}
For any $a,b,b'\in S$ and $z\in Z$
\[\tag{$Z$-Invariance}\label{ZInvariance}az=a<_{b'}b\qquad\Rightarrow\qquad a<_{b'}bz\quad\text{and}\quad a<_{zb'}b.\]
\end{prp}

\begin{proof}
If $az=a<_{b'}b$ then $b'bz,bzb'\in Z$ and $ab'bz=az$, i.e. $a<_{b'}bz$.  Also $azb'=ab'\in N$ and $azb'b=a$, i.e. $a<_{zb'}b$.
\end{proof}

We can even split up pairs whose products lie in $Z$.

\begin{prp}
For any $a,b,b',c,c'\in S$ with $cc',c'c\in Z$,
\[\tag{$Z$-Splitting}\label{ZSplitting}acc'=a<_{b'}b\qquad\Rightarrow\qquad ac<_{c'b'}bc.\]
\end{prp}

\begin{proof}
If $acc'=a<_{b'}b$ and $cc'\in Z$ then $c'b'bc,bcc'b'\in Z$, $acc'b'bc=ac$ and $acc'b'=ab'bcc'b'\in NZ\subseteq N$, i.e. $ac<_{c'b'}bc$.
\end{proof}

For use in the next section, we denote the up-closure of $A\subseteq S$ by
\[A^<=\{b>a:a\in A\}.\]

\section{Duals}\label{Duals}

\begin{dfn}
The \emph{dual} of $A\subseteq S$ is defined by
\[A^*=\{s\in S:A\ni a<_sb\}.\]
\end{dfn}

So $s\in A^*$ precisely when some $a\in A$ is dominated by another element via $s$.  In particular, if $A\subseteq B$ then $A^*\subseteq B^*$, a simple fact we will often use below.

In inverse semigroups, duals are just up-closures of inverses.

\begin{prp}
If $S$ is an inverse semigroup and $N=Z=\mathsf{E}(S)$ then
\[A^*=A^{-1<}.\]
\end{prp}

\begin{proof}
If $s\in A^*$ then we have $a\in A$ with $a<_sb$ and hence $a<_{b^{-1}}b$, by \eqref{CanonicalOrder}, so
\[aa^{-1}=ab^{-1}ba^{-1}=ab^{-1}(ab^{-1})^{-1}=ab^{-1}=asbb^{-1}=bb^{-1}as=as. \]
Then $a^{-1}=a^{-1}aa^{-1}=a^{-1}as\in Ns$ and hence $a^{-1}<s$, showing that $A^*\subseteq A^{-1<}$.

Conversely, if $s\in A^{-1<}$ then we have $a\in A$ with $a^{-1}<s$.  Again \eqref{CanonicalOrder} yields $a^{-1}<_{s^{-1}}s$.  Taking inverses yields $a<_ss^{-1}$ so $s\in A^*$, showing $A^{-1<}\subseteq A^*$.
\end{proof}

Again we want to show that duals in structured semigroups still behave much like up-closures of inverses.  For example, from \eqref{Multiplicativity}, we immediately see that the $^*$ operation is a kind of antimorphism, i.e. for all $A,B\subseteq S$.
\[\tag{Antimorphism}\label{Antimorphism}A^*B^*\subseteq(BA)^*.\]
Next we note some relations between the $^*$ and $^<$ operations.

\begin{prp}\label{Astar}
For any $A\subseteq S$,
\begin{align}
\label{<*}A^{<*}&\subseteq A^{*<}.\\
\label{**}A^{<<}&\subseteq A^{**}.
\end{align}
Consequently $A^*\neq\emptyset$ whenever $\emptyset\neq A\subseteq A^<$.
\end{prp}

\begin{proof}\
\begin{itemize}
\item[\eqref{<*}] If $c'\in A^{<*}$ then we have $a,b,b',c\in S$ with $A\ni a<_{b'}b<_{c'}c$.  Then $a<_{b'bc'}c$, by \eqref{Transitivity} and \eqref{ZInvariance}, so $A^*\ni b'bc'<_cc'\in A^{*<}$, by \eqref{Switch}.

\item[\eqref{**}] If $c\in A^{<<}$ then we have $a,b',b,c'\in S$ with $A\ni a<_{b'}b<_{c'}c$.  Then again $A^*\ni b'bc'<_cc'$ and hence $c\in A^{**}$.
\end{itemize}
For the last statement, note that if $A^*=\emptyset$ and $A\subseteq A^<$ then \eqref{**} would yield $A\subseteq A^{<<}\subseteq A^{**}\subseteq\emptyset^*=\emptyset$.
\end{proof}

If $A$ is closed under the `triple product', these inclusions become equalties.

\begin{prp}
If $AA^*A\subseteq A$ then
\begin{align}
\label{<*=}A^{<*}&=A^{*<}\subseteq A^*.\\
\label{**=}A^{<<}&=A^{**}.
\end{align}
\end{prp}

\begin{proof}\
\begin{itemize}
\item[\eqref{<*=}] Take $a'\in A^{*<}$ so we have $a,b,b',c\in S$ with $A\ni c<_{b'}b$ and $b'<_aa'$.  As $ab'ba'\in Z$, $ba'ab'=bb'\in Z$ and $bb'c=c$, \eqref{ZSplitting} and \eqref{Switch} yield
\[ab'c<_{b'ba'}ab'b<_{a'}a.\]
Then \eqref{NInvariance} yields $cb'ab'c<ab'b<_{a'}a$ and \eqref{Transitivity} yields $cb'ab'c<_{a'}a$.  Thus $a'\in A^{<*}\cap A^*$, as
\[cb'ab'c\in AA^*A^{**}A^*A\subseteq A(AA^*A)^*A\subseteq AA^*A\subseteq A,\]
showing that $A^{*<}\subseteq A^{<*}\cap A^*$ and hence $A^{<*}=A^{*<}\subseteq A^*$, by \eqref{<*}.

\item[\eqref{**=}] If we take $a\in A^{**}$ then we likewise have $a',b,b',c\in S$ with $A\ni c<_{b'}b$ and $b'<_aa'$, which again yields $cb'ab'c<ab'b<a$ and hence $a\in A^{<<}$.  Combined with \eqref{**}, this shows that $A^{<<}=A^{**}$.\qedhere
\end{itemize}
\end{proof}

\subsection{Diagonality}\label{Diagonality}

We can improve some results when $N$ is `diagonal'.

\begin{dfn}
We call $D\subseteq S$ \emph{diagonal} if, for all $a,d,b\in S$,
\[\tag{Diagonal}\label{Diagonal}ad,d,db\in D\qquad\Rightarrow\qquad adb\in D.\]
\end{dfn}

\begin{xpl}
If $S$ is a semigroup of slice-sections of a groupoid bundle $\pi:F\rightarrow G$ then the canonical diagonal $N=\{n\in S:\mathrm{dom}(n)\subseteq G^0\}$
is indeed diagonal in $S$.  To see this just note that if $an,n,nb\in N$ then $anb\in N$ because
\[\mathrm{dom}(anb)=\mathrm{dom}(a)\mathrm{dom}(n)\mathrm{dom}(n)\mathrm{dom}(b)=\mathrm{dom}(an)\mathrm{dom}(nb)\subseteq G^0.\]
\end{xpl}

\begin{xpl}
Any subsemigroup $I\subseteq S$ of idempotents (e.g. all idempotents in an inverse semigroup $S$) is diagonal, as $ai,i,ib\in I$ implies $aib=aiib\in II\subseteq I$.

More generally, any regular subsemigroup $R\subseteq S$ (i.e. for every $r\in R$ we have $r'\in R$ with $rr'r=r$) is diagonal, as $ar,r,rb\in R$ yields $r'\in R$ with
\[arb=arr'rb\in RRR\subseteq R.\]
\end{xpl}

\begin{xpl}
The range of any `conditional expectation' $\Phi$ satisfies \eqref{Diagonal} \textendash\, if $R=\mathrm{ran}(\Phi)$ for an idempotent map $\Phi$ on $S$ such that, for all $a,b\in S$,
\[\Phi(\Phi(a)b)=\Phi(a)\Phi(b)=\Phi(a\Phi(b)),\]
then $ar,r,rb\in R$ implies $\Phi(arb)=ar\Phi(b)=a\Phi(rb)=arb$ and hence $arb\in R$.
\end{xpl}

Diagonality allows us to exchange the subscript and right argument of $<$ in a slightly more general situation than in \eqref{Switch}.

\begin{prp}
If $N$ is diagonal then, for all $a,b,c,d,d'\in S$,
\[\tag{Exchange}\label{Exchange}ab,bc\in N\quad\text{and}\quad b<_{d'}d\qquad\Rightarrow\qquad abc<_dd'.\]
\end{prp}

\begin{proof}
If $b<_{d'}d$ and $ab,bc\in N$ then $abcdd'=add'bc=abc$ and \eqref{Diagonal} yields $abcd=add'bcd\in N$, as $ab=add'b,d'b,d'bcd\in N$, i.e. $abc<_dd'$.
\end{proof}

This yields is a kind of transitivity, namely
\[\tag{*-Transitivity}\label{*Transitivity}a<_{b'}b\quad\text{and}\quad b'<_cc'\qquad\Rightarrow\qquad a<_{c'}c.\]
Indeed, if $a<_{b'}b$ and $b'<_cc'$ then \eqref{Exchange} yields $a=ab'b<_{c'}c$.

We can then prove \eqref{<*=} and \eqref{**=} even without the triple product assumption.

\begin{prp}
If $N$ is diagonal then, for any $C\subseteq S$,
\begin{equation}\label{DiagonalStars}
C^{*<}\subseteq C^*,\qquad C^{**}=C^{<<}\qquad\text{and}\qquad C^{***}\subseteq C^*.
\end{equation}
\end{prp}

\begin{proof}
If $a'\in C^{*<}$ then we have $a,b,b',c\in S$ with $C\ni c<_{b'}b$ and $b'<_aa'$.  By \eqref{*Transitivity}, $c<_{a'}a$ and hence $a'\in C^*$, showing $C^{*<}\subseteq C^*$.

If $a\in C^{**}$ then we have $a',b,b',c\in S$ with $C\ni c<_{b'}b$ and $b'<_aa'$.  By \eqref{Exchange}, $c<_{b'}bb'b<_{a'}a$ and hence $a\in C^{<<}$, showing $C^{**}\subseteq C^{<<}$.  The reverse inclusion was already proved in \eqref{**}.

Now it follows that $C^{***}\subseteq C^{*<<}\subseteq C^{*<}\subseteq C^*$.
\end{proof}

\section{Atlases}\label{Atlases}

\begin{dfn}
We call $A\subseteq S$ an \emph{atlas} and $C\subseteq S$ a \emph{coset} if
\begin{align}
\tag{Atlas}AA^*A\subseteq A&\subseteq A^<.\\
\tag{Coset}CC^*C\subseteq C&=C^<.
\end{align}
\end{dfn}

So an atlas is a triple-product-closed round subset, and a coset is also an up-set.  These generalise the same notions for inverse semigroups from \cite[\S1.4 before Proposition 26]{Lawson1998} (as the order there is reflexive, all subsets are trivially round).

\begin{rmk}\label{PrincipalFilters}
When $S$ is a group with identity $e$ and $N=Z=\{e\}$, a coset $C$ of $S$ is precisely a coset in the usual sense, specifically a left coset of the subgroup $C^{-1}C$ or a right coset of the subgroup $CC^{-1}$.  In particular, each singleton $\{a\}$ is a coset of the trivial subgroup $\{e\}$.

More generally, if $S$ is an inverse semigroup and $N=Z=\mathsf{E}(S)$ then every principal filter $a^\leq$ is a coset.  If $S$ is the normaliser semigroup of a Cartan subalgebra $C$ of some C*-algebra and $N=Z=C$ then again we have principal filters $a^<$, for each $a\in S$, which are also cosets, by \autoref{FilterCosets} below.  We also have principal filters in the more general structured C*-algebras considered in \cite{Bice2021DHKR}, thanks to \cite[Proposition 6.11]{Bice2021DHKR} and \cite[Lemmas 5.2 and 5.3]{Bice2021Pierce}.

However, general structured semigroups can have few cosets, e.g. $S$ itself is the only non-empty coset in \autoref{NonFaithfulExample} below.  Indeed, an important consequence of one of our results is that we have faithful bundle representations precisely when there are enough cosets to distinguish the elements of $S$ via their corresponding equivalence relations, at least when $Z$ is symmetric -- see \autoref{CosetBundleFaithful} below.
\end{rmk}

First we show that atlases generate cosets.

\begin{prp}\label{Atlas->Cosets}
If $A$ is an atlas then $A^*$ and $A^<=A^{**}(\supseteq A)$ are cosets.
\end{prp}

\begin{proof}
If $A$ is an atlas, \eqref{Antimorphism} yields $A^*A^{**}A^*\subseteq(AA^*A)^*\subseteq A^*$.  Also
\[A^*\subseteq A^{<*}=A^{*<}\subseteq A^*,\]
by \eqref{<*=}, i.e. $A^*=A^{*<}$.  This shows that $A^*$ is a coset and hence $A^{**}=A^{<<}=A^<$ is also a coset, by \eqref{Transitivity}, $A\subseteq A^<$ and \eqref{**=}.
\end{proof}

\begin{prp}
If $A,B\subseteq S$ are atlases and $c\in S$ then
\begin{equation}\label{AcBc}
A\subseteq B\quad\text{and}\quad A^*=B^*\qquad\Rightarrow\qquad(Ac)^*=(Bc)^*.
\end{equation}
\end{prp}

\begin{proof}
Assume $A,B\subseteq S$ are atlases with $A\subseteq B$ and $A^*=B^*$.  If $d\in(Bc)^*$, then we have $b\in B$ and $d'\in S$ with $bc<_dd'$.  As $A\subseteq A^<$ and $B\subseteq B^<$, we can take $n\in AA^*\cap N$, $z\in A^*A\cap N$ and $b'\in B^*$ with $b'b\in Z$.  As $Z$ is binormal, $nbzb'\in NZ\subseteq N$ and hence \eqref{NInvariance} yields $nbzb'bc<_dd'$.  Also
\[nbzb'b=nbb'bz\in AA^*BB^*BA^*A=AB^*BB^*BB^*A\subseteq AB^*A=AA^*A\subseteq A,\]
so $d\in(Ac)^*$, showing that $(Bc)^*\subseteq(Ac)^*\subseteq(Bc)^*$.
\end{proof}

\begin{dfn}\label{ActionNotation}
For any $A\subseteq S$ and $b\in S$ we define
\begin{align*}
A|b\qquad&\Leftrightarrow\qquad\exists b'\in S\ (bb',b'b\in Z\text{ and }\exists a\in A\ (a=abb')).\\
b|A\qquad&\Leftrightarrow\qquad\exists b'\in S\ (bb',b'b\in Z\text{ and }\exists a\in A\ (a=b'ba)).
\end{align*}
\end{dfn}

We read $A|b$ and $b|A$ as saying $b$ \emph{acts on} $A$ (from the right or left respectively).  Indeed, when $b$ acts on $A$, the product is another atlas \textendash\, see \autoref{A|b} below.

\begin{prp}\label{ABatlas}
If $A,B\subseteq S$ are atlases satisfying $A^*A\subseteq(BB^*)^<$ then their product $AB$ is also an atlas.  Consequently $(AB)^*$ is a coset and, for all $a\in A$,
\begin{equation}\label{BC=bC}
a|B\qquad\text{and}\qquad(aB)^*=(AB)^*.
\end{equation}
\end{prp}

\begin{proof}
If $A^*A\subseteq(BB^*)^<$ then \eqref{Multiplicativity} and the fact $B$ is an atlas yields $A^*AB\subseteq(BB^*)^<B^<\subseteq(BB^*B)^<\subseteq B^<$ and hence
\[(AB)^*A\subseteq(AB)^*A^{**}\subseteq(A^*AB)^*\subseteq B^{<*}=B^{*<}\subseteq B^*.\]
Thus $AB(AB)^*AB\subseteq ABB^*B\subseteq AB\subseteq A^<B^<\subseteq(AB)^<$, by \eqref{Multiplicativity}, so $AB$ is an atlas.

Take $a',c\in S$ with $A\ni c<_{a'}a$.  Then $a'c\in A^*A\subseteq(BB^*)^<$ so, taking any $b\in B$, $a'cb\in(BB^*)^<B\subseteq B^<$.  Taking $e,d\in S$ with $B\ni e<_da'cb$, we see that
$a'ae=a'aa'cbde=a'cbde=e$ so $a|B$.  For any other $f\in A$, \eqref{NInvariance} yields
\[(fB)^*\subseteq(aa'fB)^*\subseteq(aA^*AB)^*\subseteq(aB^<)^*=(aB)^*,\]
by \eqref{AcBc} (applied to $S^\mathrm{op}$ with $a$, $B^<$ and $B$ in place of $c$, $A$ and $B$ respectively).  This shows that $(AB)^*\subseteq(aB)^*\subseteq(AB)^*$.
\end{proof}

In particular, we can take $B=A^*$ above, as $A^*A\subseteq A^{<*}A^<\subseteq(A^*A)^<$, i.e.
\[A\text{ is an atlas}\qquad\Rightarrow\qquad AA^*\text{ is an atlas}.\]
In this case, the cosets $(AA^*)^*$ and $(AA^*)^<$ they generate are the same as $(AA^*)^*\subseteq(A^{**}A^*)^*\subseteq(AA^*)^{**}=(AA^*)^<$ and $(AA^*)^<\subseteq(A^{**}A^*)^<\subseteq(AA^*)^{*<}=(AA^*)^*$.  Likewise $(A^*A)^*$ and $(A^*A)^<$ coincide, and we denote these cosets by
\[\mathsf{s}(A)=(A^*A)^*=(A^*A)^<\qquad\text{and}\qquad\mathsf{r}(A)=(AA^*)^*=(AA^*)^<.\]
We will soon see that these are indeed the source and range in a groupoid of cosets.

\begin{prp}\label{A|b}
If $A|b$ and $A$ is an atlas then $Ab$ is an atlas with $\mathsf{r}(Ab)=\mathsf{r}(A)$.
\end{prp}

\begin{proof}
Take $b'\in S$ with $bb',b'b\in Z$ and $a\in A$ with $a=abb'$.  As $A\subseteq A^<$, we can also take $a'\in A^*$ with $aa',a'a\in Z$.

First note
\begin{equation}\label{b'A*}
b'A^*\subseteq(Ab)^*.
\end{equation}
To see this, take any $c'\in A^*$, so we have $c,d\in S$ with $A\ni d<_{c'}c$.  Then \eqref{NInvariance} yields $da'a<_{c'}c$ and \eqref{ZSplitting} yields $da'ab<_{b'c'}cb$.  Thus $b'c'\in(Ab)^*$, as $da'a\in AA^*A\subseteq A$, showing that $b'A^*\subseteq(Ab)^*$.

Next we claim that
\begin{equation}
b(Ab)^*\subseteq A^*.
\end{equation}
To see this, take $c'\in(Ab)^*$, so we have $c\in S$ and $d\in A$ with $db<_{c'}c$.  Take $d'\in A^*$ with $dd',d'd\in Z$ and let $n=da'ad'\in N$ so \eqref{NInvariance} yields $ndb<_{c'}c$.  Also $ndbb'=da'ad'dbb'=da'abb'd'd=da'ad'd=nd$ and hence $ndbb'b=ndb$ so \eqref{ZSplitting} yields $nd=ndbb'<_{bc'}cb'$.  Thus $bc'\in A^*$, as $nd\in AA^*AA^*A\subseteq A$, showing that $b(Ab)^*\subseteq A^*$ and hence
\begin{equation}\label{Ab3x}
Ab(Ab)^*Ab\subseteq AA^*Ab\subseteq Ab.
\end{equation}

Next note $a'aa'=a'abb'a'=bb'a'aa'$, so $b$ and $a'aa'\in A^*AA^*\subseteq A^*$ witness $b'|A^*$.  Thus \eqref{b'A*} applied in $S^\mathrm{op}$ yields $A^{**}b\subseteq(b'A^*)^*$.  Then \eqref{**=} and \eqref{Ab3x} yield
\[Ab\subseteq A^{<<}b=A^{**}b\subseteq(b'A^*)^*\subseteq(Ab)^{**}=(Ab)^{<<}\subseteq(Ab)^<,\]
showing that $Ab$ is an atlas.

Finally, note $Ab(Ab)^*\subseteq AA^*\subseteq A^{**}(Abb')^*\subseteq(Abb'A^*)^*\subseteq(Ab(Ab)^*)^*=\mathsf{r}(Ab)$ so taking duals yields $\mathsf{r}(Ab)\subseteq\mathsf{r}(A)\subseteq\mathsf{r}(Ab)^*=\mathsf{r}(Ab)$.
\end{proof}

\begin{prp}
If $A\subseteq S$ is an atlas and $n\in N$ then
\begin{equation}\label{A|e}
A|n\qquad\Leftrightarrow\qquad n\in\mathsf{s}(A)\qquad\Rightarrow\qquad(An)^<=A^<.
\end{equation}
\end{prp}

\begin{proof}
If $n\in\mathsf{s}(A)$ then $An\subseteq A^<(A^*A)^<\subseteq(AA^*A)^<\subseteq A^<$ and hence
\[A^<\subseteq(An)^<\subseteq A^{<<}\subseteq A^<.\]
Also, we have $m,n'\in S$ with $A^*A\ni m<_nn'$ so $nn',n'n\in Z$ and $mnn'=m$.  Taking any $a\in A$, note $am\in AA^*A\subseteq A$ and $amnn'=am$, showing that $A|n$.  

Conversely, if $A|n$, we have $n'\in S$ with $nn',n'n\in Z$ and $a\in A$ with $a=ann'$.  Then $a'a=a'ann'$, for any $a'\in A^*$ with $a'a\in N$, so $a'a<_nn'$ and $n\in\mathsf{s}(A)$.
\end{proof}

\section{Cosets}\label{Cosets}

Our goal here is to show that the non-empty cosets
\[\mathcal{C}(S)=\{C\subseteq S:CC^*C\subseteq C=C^<\neq\emptyset\}\]
form an \'etale groupoid.  This generalises similar results for filters in inverse semigroups in \cite{Lenz2008} and \cite{LawsonMargolisSteinberg2013} (we will have more to say about filters in \autoref{Filters}).

\begin{thm}
$\mathcal{C}(S)$ is a groupoid under the inverse $C\mapsto C^*$ and product
\[B\cdot C=(BC)^<\quad\text{when}\quad\mathsf{s}(B)=\mathsf{r}(C).\]
\end{thm}

\begin{proof}
If $(B,C)\in\mathcal{C}^2=\{(B,C):B,C\in\mathcal{C}(S)\ \text{and}\ \mathsf{s}(B)=\mathsf{r}(C)\}$ then $(BC)^<$ is coset, by \autoref{Atlas->Cosets} and \autoref{ABatlas}.  Also, as $B\neq\emptyset\neq C$, it follows that $\emptyset\neq BC\subseteq B^<C^<\subseteq(BC)^<$, so the product is well-defined on $\mathcal{C}(S)$.  Moreover, if $c\in C$, \autoref{ABatlas} yields $B|c$ and $(BC)^<=(Bc)^<$ so \autoref{A|b} yields
\[\mathsf{r}((BC)^<)=\mathsf{r}((Bc)^<)=\mathsf{r}(Bc)=\mathsf{r}(B).\]
Thus $(A,B)\in\mathcal{C}^2$ iff $(A,(BC)^<)\in\mathcal{C}^2$ in which case, for any $c\in C$,
\[(A(BC)^<)^<=(Abc)^<\subseteq(ABC)^<\subseteq(A(BC)^<)^<.\]
Likewise, for any $(A,B)\in\mathcal{C}^2$, we see that $\mathsf{s}((AB)^<)=\mathsf{s}(B)$ so $((AB)^<,C)\in\mathcal{C}^2$ iff $(B,C)\in\mathcal{C}^2$ in which case $((AB)^<C)^<=(ABC)^<=(A(BC)^<)^<$, showing that the product is associative.

Again by \autoref{Atlas->Cosets}, as well as the last part of \autoref{Astar}, if $C\in\mathcal{C}(S)$ then $C^*\in\mathcal{C}(S)$ so the involution is well-defined on $\mathcal{C}(S)$.  Also, for any $(B,C)\in\mathcal{C}^2$, $c\in C$ and $c'\in C^*$ with $cc'\in N$,
\[B\subseteq(Bcc')^<=(BCC^*)^<\subseteq(B(B^*B)^<)^<\subseteq B,\]
showing that $C\cdot C^*$ is a unit.  Likewise, $(B^*BC)^<=C$ so $B^*\cdot B$ is a unit, showing that the involution takes each element to its inverse.  Thus $\mathcal{C}(S)$ is a groupoid.
\end{proof}

The units of this groupoid have a couple of simple characterisations.

\begin{prp}\label{Cunit}
$C\in\mathcal{C}(S)$ is a unit iff $C\cap N\neq\emptyset$ iff $C\cap Z\neq\emptyset$.
\end{prp}

\begin{proof}
If $C\in\mathcal{C}(S)$ is a unit then $C=(CC^*)^<\supseteq CC^*$.  Taking any $b,b',c\in S$ with $C\ni c<_{b'}b$, we see that $bb'\in CC^*\cap Z\subseteq C\cap Z$.  Conversely, if $n\in C\cap N$ then, for any $B$ with $(B,C)\in\mathcal{C}^2$, \autoref{ABatlas} and \eqref{A|e} yield $B=(Bn)^<=(BC)^<$.  Likewise, $B=(CB)^<$ for any $B$ with $(C,B)\in\mathcal{C}^2$ so $C$ is a unit in $\mathcal{C}(S)$.
\end{proof}

The following slices of $\mathcal{C}(S)$ will play an important role very soon.

\begin{prp}\label{CaSlice}
For every $a\in S$, we have a slice in $\mathcal{C}(S)$ given by
\[\mathcal{C}_a=\{C\in\mathcal{C}(S):a\in C\}.\]
\end{prp}

\begin{proof}
If $B,C\in\mathcal{C}_a$ and $\mathsf{s}(B)=\mathsf{s}(C)$ then, by \eqref{BC=bC},
\[B=B\cdot\mathsf{s}(B)=(a\mathsf{s}(B))^<=(a\mathsf{s}(C))^<=C\cdot\mathsf{s}(C)=C.\]
Likewise, $\mathsf{r}(B)=\mathsf{r}(C)$ implies $B=C$ so $\mathcal{C}_a$ is a slice.
\end{proof}

The canonical topology on $\mathcal{C}(S)$ is generated by the slices $(\mathcal{C}_a)_{a\in S}$.  So a basis for this topology is given by $\mathcal{C}_F=\bigcap_{f\in F}\mathcal{C}_f$, for finite $F\subseteq S$.

\begin{thm}\label{EtaleCosets}
The coset groupoid $\mathcal{C}(S)$ is \'etale in the canonical topology.
\end{thm}

\begin{proof}
If $C^*\in\mathcal{C}_b$ then $b\in C^*$ so we have $a,b'\in S$ with $C\ni a<_bb'$ and hence $C\in\mathcal{C}_a$ and $\mathcal{C}_a^*\subseteq\mathcal{C}_b$, showing that the involution $C\mapsto C^*$ is continuous on $\mathcal{C}(S)$.  Similarly, if $B\cdot C\in\mathcal{C}_a$ then $a\in B\cdot C$ so we have $b\in B$ and $c\in C$ with $bc<a$, i.e. $B\in\mathcal{C}_b$, $C\in\mathcal{C}_c$ and $\mathcal{C}_b\cdot\mathcal{C}_c\subseteq\mathcal{C}_a$, showing that the product is also continuous.

To see that the source $\mathsf{s}$ is an open map on $\mathcal{C}(S)$, take any finite $F\subseteq S$ and $C\in\mathcal{C}_F$.  Further take finite $G\subseteq C$ with $F\subseteq G^<$ and fix some $a,b',b\in S$ with $G\ni a<_{b'}b\in F$.  We claim that
\[\mathsf{s}(C)\in\mathcal{C}_{b'G}\subseteq\mathsf{s}[\mathcal{C}_F].\]
To see this, first note that $b'\in G^*\subseteq C^*$ and $G\subseteq C$ implies $b'G\subseteq C^*C\subseteq\mathsf{s}(C)$ and hence $\mathsf{s}(C)\in\mathcal{C}_{b'G}$.  Next note that if $B\in\mathcal{C}_{b'G}$ then $b'a\in B\cap N$ so $B$ is a unit, by \autoref{Cunit}, and $b'bb'a=b'a$ so $b|B$ and $(bB)^<\in\mathcal{C}_{bb'G}\subseteq\mathcal{C}_F$ (as $A\in\mathcal{C}_{bb'G}$ implies $bb'G\subseteq A$ and hence $F\subseteq G^<\subseteq(bb'G)^<\subseteq A^<\subseteq A$).  Thus $B=\mathsf{s}(B)=\mathsf{s}(bB)\in\mathsf{s}[\mathcal{C}_F]$, showing that $\mathcal{C}_{b'G}\subseteq\mathsf{s}[\mathcal{C}_F]$.  As $C$ was arbitrary, this shows $\mathsf{s}[\mathcal{C}_F]$ is open which, as $F$ was arbitrary, shows that $\mathsf{s}$ is an open map.
\end{proof}

\subsection{Symmetry}\label{Symmetry}

The structured semigroups we are interested in often satisfy a certain additional assumption, which will also allow us to say more about cosets.

\begin{dfn}
We call $Y\subseteq S$ \emph{symmetric} if, for all $a,b\in S$,
\[\tag{Symmetric}\label{Symmetric}ab\in Y\qquad\Rightarrow\qquad baba\in Y.\]
\end{dfn}

For example, the idempotents $\mathsf{E}(S)$ in any semigroup $S$ are always symmetric, as $ab\in\mathsf{E}(S)$ implies that $b(ababab)a=baba\in\mathsf{E}(S)$.  Likewise, one can verify that the slice-sections $\mathcal{S}(\pi|_{F^\times})$ of the invertible part of a Fell bundle $\pi$ have a symmetric central-diagonal $\mathcal{Z}(\pi|_{F^\times})$ (see \autoref{FellBundles}).

Symmetry yields even more one-sided versions of $<$ than in \autoref{1sided<}.

\begin{prp}\label{very1sided<}
If $Z$ is symmetric and $a,b\in S$ then
\begin{align}
asb=a,\ as\in N\text{ and }(bs\in Z\text{ or }sb\in Z)\quad&\Rightarrow\quad a<_sbsb\text{ and }a<_{sbs}b.\\
bsa=a,\ sa\in N\text{ and }(bs\in Z\text{ or }sb\in Z)\quad&\Rightarrow\quad a<_sbsb\text{ and }a<_{sbs}b.
\end{align}
\end{prp}

\begin{proof}
If $bs\in Z$ then $bsbs\in ZZ\subseteq Z$ and $sbsb\in Z$, by symmetry.  If $sb\in Z$ then, likewise, $sbsb\in ZZ\subseteq Z$ and $bsbs\in Z$, by symmetry.  If $asb=a$ and $as\in N$ too then $asbs\in NZ\subseteq N$ and $asbsb=asb=a$ so $a<_sbsb$ and $a<_{sbs}b$.  This proves the first $\Rightarrow$ and the second follows by a dual argument.
\end{proof}

\begin{prp}
If $Z$ is symmetric then, for all $a,b\in S$ and $n\in N$,
\begin{align}
\label{a<be}a<_{b'}nb\qquad&\Rightarrow\qquad a<_{b'nbb'n}b.\\
\intertext{If $N$ is diagonal and $Z$ is symmetric then, for all $a,b\in S$ and $n\in N$,}
\label{a<nb}a<_{b'n}b\qquad&\Rightarrow\qquad a<_{b'}nbb'nb.
\end{align}
\end{prp}

\begin{proof}
If $a<_{b'}nb$ then $ab'n\in NN\subseteq N$, $ab'nb=a$ and $b'nb\in Z$ and hence $a<_{b'nbb'n}b$, by \autoref{very1sided<} with $s=b'n$.

Now say $N$ is also diagonal and $a<_{b'n}b$.  In particular, $bb'n\in Z$ so, by symmetry, $nbb'nbb'\in Z\subseteq N$.  As $ab'n,n\in N$ too, diagonality yields $ab'=ab'nbb'nbb'\in N$.  As $ab'nb=a$ and $b'nb\in Z$, we then see that $a<_{b'}nbb'nb$, again by \autoref{very1sided<} but with $s=b'$ and $nb$ replacing $b$.
\end{proof}

\begin{prp}
If $Z$ is symmetric then, for any atlas $A\subseteq S$ and $a,b',b\in S$,
\begin{equation}\label{SymmetricAction}
abb'=a\in A\quad\text{and}\quad bb'\in Z\qquad\Rightarrow\qquad A|b.
\end{equation}
\end{prp}

\begin{proof}
If $bb'\in Z$ then $b'bb'b\in Z$ by symmetry.  Also $bb'bb'\in ZZ\subseteq Z$ and $abb'bb'=abb'=a$, showing that $b'bb'$ witnesses $A|b$.
\end{proof}

Let us now extend our notion of \'etale representation from \autoref{InverseSemigroups}.

\begin{dfn}\label{EtaleRep}
An \emph{\'etale representation} of the structured semigroup $(S,N,Z)$ on an \'etale groupoid $G$ is a semigroup homomorphism $\theta:S\rightarrow\mathcal{B}(G)$ to the open slices of $G$ such that $\theta[S]$ covers $G$, $\theta[N]\subseteq\mathcal{O}(G^0)$ and, for all $g\in G$ and $a\in S$,
\begin{equation}
\label{RoundRep}\tag{Locally Round}g\in\theta(a)\qquad\Rightarrow\qquad\exists b<a\ (g\in\theta(b)).
\end{equation}
\end{dfn}

Intuitively, \eqref{RoundRep} is saying we can always shrink open neighbourhoods arising from the representation.  On a more technical level, we must restrict the representations in this way for the coset representation to be universal, as shown below.  The following observation shows that we can at least rest assured that requiring \eqref{RoundRep} and $\theta[N]\subseteq\mathcal{O}(G^0)$ is still consistent with the original notion of an \'etale representation of an inverse semigroup introduced in \autoref{InverseSemigroups}.

\begin{prp}
If $S$ is an inverse semigroup and $N=Z=\mathsf{E}(S)$ then any semigroup homomorphism $\theta:S\rightarrow\mathcal{B}(G)$ is automatically an \'etale representation.
\end{prp}

\begin{proof}
As $\theta$ is semigroup homomorphism, it maps idempotents in $S$ to idempotents in $\mathcal{B}(G)$, which are precisely the open subsets of the unit space $G^0$, i.e. $\theta[\mathsf{E}(S)]\subseteq\mathcal{O}(G^0)$.  Moreover, the domination relation $<$ here is just the usual order on the inverse semigroup $S$ -- see \eqref{CanonicalOrder}.  In particular, $<$ is reflexive and hence \eqref{RoundRep} is immediately verified by taking $b=a$.
\end{proof}

As before in \autoref{InverseSemigroups}, we call an \'etale representation $\mu:S\rightarrow\mathcal{B}(G)$ \emph{universal} if, for every \'etale representation $\theta:S\rightarrow\mathcal{B}(H)$, there exists a unique \'etale morphism $\phi:H\rightarrow G$ such that $\theta=\overline\phi\circ\mu$, where $\overline\phi:\mathcal{B}(G)\rightarrow\mathcal{B}(H)$ is the preimage map.

Again as before, we let $\mathcal{C}$ denote the map $a\mapsto\mathcal{C}_a$.

\begin{thm}\label{SymmetricEtaleRepresentation}
If $Z$ is symmetric then $\mathcal{C}$ is a universal \'etale representation.
\end{thm}

\begin{proof}
First we must show that $\mathcal{C}$ is a semigroup homomorphism.  For $a,b\in S$, certainly $\mathcal{C}_a\cdot\mathcal{C}_b\subseteq\mathcal{C}_{ab}$.  Next, for $n\in N$, we have a partial converse, namely
\begin{equation}\label{CaeCa}
\mathcal{C}_{an}\subseteq\mathcal{C}_a.
\end{equation}
Indeed, if $C\in\mathcal{C}_{an}$ then we have $c\in C$ with $c<an$.  This implies $c<a$, by \eqref{a<be}, and hence $a\in C^<\subseteq C$, i.e. $C\in\mathcal{C}_a$, which proves \eqref{CaeCa}.

More generally, if $C\in\mathcal{C}_{ab}$ then we have $c,c'\in S$ with $C\ni c<_{c'}ab$.  In particular, $abc'\in Z$ and $abc'c=c$ and hence $bc'|C$, by \eqref{SymmetricAction}.  By \autoref{A|b}, $B=(bc'C)^<$ is a coset with $\mathsf{s}(B)=\mathsf{s}(C)$.  Also $B\in\mathcal{C}_{bc'c}\subseteq\mathcal{C}_b$, by \eqref{CaeCa}, as $c'c\in N$.  As $b\in B$, we have $b'\in B^*$ with $bb'\in Z\subseteq N$ and hence $C\cdot B^*\in\mathcal{C}_{ab}\cdot\mathcal{C}_{b'}\subseteq\mathcal{C}_{abb'}\subseteq\mathcal{C}_a$, again by \eqref{CaeCa}.  Thus $C=C\cdot B^*\cdot B\in\mathcal{C}_a\cdot\mathcal{C}_b$.  This shows that $\mathcal{C}_{ab}=\mathcal{C}_a\cdot\mathcal{C}_b$, for all $a,b\in S$, so $\mathcal{C}$ is indeed a semigroup homomorphism.

To see that $\mathcal{C}$ is even an \'etale representation, note first that $n\in N$ implies that $\mathcal{C}_n$ consists entirely of unit cosets, thanks to \autoref{Cunit}.  Also any coset $C$ must satisfy $C\subseteq C^<$ by definition, so if $C\in\mathcal{C}_a$ then we have $b\in C$ with $b<a$ and hence $C\in\mathcal{C}_b$, i.e. \eqref{RoundRep} is also satisfied.  As each $C\in\mathcal{C}(S)$ is non-empty, we have $a\in C$ and hence $C\in\mathcal{C}_a$, i.e. $(\mathcal{C}_a)_{a\in S}$ covers $\mathcal{C}(S)$, as required.

To see that $\mathcal{C}$ is universal, take any other \'etale representation $\theta:S\rightarrow\mathcal{B}(G)$.  For each $g\in G$, define $\phi(g)\subseteq S$ by
\[\phi(g)=\{a\in S:g\in\theta(a)\}.\]
We claim $\phi(g)\in\mathcal{C}(S)$.  Certainly $\phi(g)\neq\emptyset$, as $\theta[S]$ covers $G$.  Also $\phi(g)\subseteq\phi(g)^<$, as $\theta$ is required to satisfy \eqref{RoundRep}.  Conversely, say $\phi(g)\ni a<b$ so $g\in\theta(a)=\theta(asb)=\theta(a)\theta(s)\theta(b)$ and $\theta(a)\theta(s)=\theta(as)\in\theta[N]\subseteq\mathcal{O}(G^0)$.  Arguing as in the proof of \autoref{<Equiv} (or considering the special case of \autoref{<Equiv} where $\pi:G\rightarrow G$ is the identity bundle and hence $\mathcal{B}(G)\approx\mathcal{S}(\pi)$), it follows that $g\in\theta(b)$ and hence $b\in\phi(g)$, showing that $\phi(g)^<\subseteq\phi(g)$.  Finally note that if $b\in\phi(g)^*$ then $a<_bc$, for some $a\in\phi(g)$, which means $g\in\theta(a)$ and hence $g^{-1}\in\theta(b)$, again by \autoref{<Equiv}.  This shows that $\phi(g)^*\subseteq\phi(g^{-1})$ so, for any $a,c\in\phi(g)$ and $b\in\phi(g)^*$, we see that $g=gg^{-1}g\in\theta(a)\theta(b)\theta(c)=\theta(abc)$, i.e. $abc\in\phi(g)$.  Thus we have shown that $\phi(g)$ is a non-empty coset, proving the claim.

Now $\phi(g)^*\subseteq \phi(g^{-1})=\phi(g^{-1})^{**}\subseteq \phi(g)^*$, i.e. $\phi(g)^*=\phi(g^{-1})$.  Also
\[\mathsf{r}(\phi(g))=\phi(g)\cdot \phi(g)^*=(\phi(g)\phi(g)^*)^<=(\phi(g)\phi(g^{-1}))^<\subseteq\phi(gg^{-1})^<=\phi(\mathsf{r}(g)).\]
As $\phi(\mathsf{r}(g))$ contains a unit coset and hence element of $N$, $\phi(\mathsf{r}(g))$ is also a unit coset.  Taking any $a\in\phi(g)$ and $a'\in\phi(g^{-1})$, we see that $\mathsf{r}(g)=gg^{-1}\in\theta(a)\theta(a')=\theta(aa')$.  Thus $aa'\in\phi(\mathsf{r}(g))$ and then \eqref{**=} and \eqref{BC=bC} yield
\[\phi(\mathsf{r}(g))=(\phi(\mathsf{r}(g))aa')^<\subseteq(\phi(\mathsf{r}(g))\phi(g)\phi(g^{-1}))^<\subseteq(\phi(\mathsf{r}(g)g)\phi(g^{-1}))^<=\mathsf{r}(\phi(g)).\]
Likewise, $\mathsf{s}(\phi(g))=\phi(\mathsf{s}(g))$.  If we have another $h\in G$ with $\mathsf{s}(g)=\mathsf{r}(h)$ then $\mathsf{s}(\phi(g))=\phi(\mathsf{s}(g))=\phi(\mathsf{r}(h))=\mathsf{r}(\phi(h))$ so $\phi(g)\cdot\phi(h)$ is defined and
\[\phi(g)\cdot\phi(h)=(\phi(g)\phi(h))^<\subseteq\phi(gh)^<=\phi(gh).\]
Likewise, $\phi(gh)\cdot\phi(h^{-1})\subseteq\phi(ghh^{-1})=\phi(g)$ and hence
\[\phi(gh)=\phi(gh)\cdot\phi(h^{-1})\cdot\phi(h)\subseteq\phi(g)\cdot\phi(h).\]
This shows that $\phi(gh)=\phi(g)\cdot\phi(h)$ and hence $\phi$ is a functor.

Now say $x\in G^0$, $C\in\mathcal{C}(S)$ and $\mathsf{r}(C)=\phi(x)$.  Taking any $c\in C$ and $c'\in C^*$ we see that $cc'\in CC^*\subseteq\mathsf{r}(C)=\phi(x)$ so $x\in\theta(cc')=\theta(c)\theta(c')$ and hence $x=\mathsf{r}(g)$, for some $g\in\theta(c)\cap\theta(c')^{-1}$.  Then \eqref{**=} and \eqref{BC=bC} again yield
\[C=\mathsf{r}(C)\cdot C=(\mathsf{r}(C)c)^<=(\phi(x)c)^<=(\phi(x)\phi(g))^<=\phi(g).\]
This shows that $\phi$ is star-surjective.  To see that $\phi$ is also star-injective, just note that $\mathsf{r}(g)=\mathsf{r}(h)$ and $\phi(g)=\phi(h)$, or even just $\phi(g)\cap\phi(h)\neq\emptyset$, implies $g=h$, as $a\in \phi(g)\cap\phi(h)$ implies that $g$ and $h$ are in the same slice $\theta(a)$.

Now note that, for any $a\in S$, $\phi^{-1}[\mathcal{C}_a]=\theta(a)\in\mathcal{B}(G)$ because
\[\phi(g)\in\mathcal{C}_a\qquad\Leftrightarrow\qquad a\in \phi(g)\qquad\Leftrightarrow\qquad g\in\theta(a).\]
Thus $\phi$ is an \'etale morphism from $G$ to $\mathcal{C}(S)$ with $\theta=\overline\phi\circ\mathcal{C}$.

Finally, for uniqueness, say we had another map $\psi:G\rightarrow\mathcal{C}(S)$ with $\overline\psi\circ\mathcal{C}=\theta$, i.e. $\psi^{-1}[\mathcal{C}_a]=\theta(a)$, for all $a\in S$.  Then again we see that $g\in\theta(a)$ iff $\psi(g)\in\mathcal{C}_a$ iff $a\in\psi(g)$ and hence $\psi(g)=\{a\in S:g\in\theta(a)\}=\phi(g)$, for all $g\in G$.
\end{proof}

\begin{rmk}
If the subbasic slices of a groupoid are closed under pointwise products then the same is true of arbitrary open slices \textendash\, see the proof of \cite[Proposition 3.18]{BiceStarling2020HTight}.  In particular, this applies to $\mathcal{C}(S)$ if $\mathcal{C}$ is a semigroup homomorphism, as shown in the first part of the proof above.  As $C\mapsto C^*$ is immediately seen to be continuous on $\mathcal{C}(S)$, as mentioned at the start of the proof of \autoref{EtaleCosets}, the open slices of $\mathcal{C}(S)$ are also closed under pointwise inverses.  It follows that $\mathcal{C}(S)$ is \'etale, again by \cite[Proposition 3.18]{BiceStarling2020HTight}, which provides an alternative proof of \autoref{EtaleCosets} in the symmetric case.
\end{rmk}

While the coset representation $\mathcal{C}$ is universal, it can only be faithful when $N$ consists entirely of idempotents.  If we want faithful representations of more general structured semigroups, including those with few idempotents, then we need to consider finer representations on groupoid bundles.  In particular, to construct such a bundle over the coset groupoid, we will split up the cosets into equivalences classes with respect to certain relations that we now proceed to examine.

\section{Equivalence}\label{Equivalence}

\begin{dfn}
For any atlas $A\subseteq S$, we define a relation $\sim_A$ on the coset $A^<$ by
\[a\sim_Ab\qquad\Leftrightarrow\qquad\exists s,t\in A^*\ (sat=sbt).\]
\end{dfn}

We can think of $a\sim_Ab$ as saying $a$ and $b$ have the same `germ' at $A$.

\begin{xpl}
Say $S=N=Z=C_0(X)=$ continuous functions from $X$ to $\mathbb{C}$ vanishing at infinity, where $X$ is a locally compact Hausdorff space.  Then each $x\in X$ defines a coset (in fact an ultrafilter)
\[S_x=\{a\in S:a(x)\neq0\}.\]
(note this is consistent with our notation from \autoref{local} when we identify each $a\in S$ with the restriction to its support $a|_{\mathrm{supp}(a)}$).  Moreover, for all $a,b\in S_x$,
\[a\sim_{S_x}b\qquad\Leftrightarrow\qquad a(y)=b(y)\text{ for all $y$ in some neighbourhood $N$ of }x.\]
Indeed, if $a\sim_{S_x}b$, i.e. $sat=sbt$, for some $s,t\in S_x^*=S_x$, then $a(x)=b(x)$, for all $x\in\mathrm{supp}(s)\cap\mathrm{supp}(t)$.  Conversely, if $a(y)=b(y)$ for all $y$ in a neighbourhood $N$ of $x$ then, taking any $s\in S_x$ with $\mathrm{supp}(s)\subseteq N$, we see that $sas=sbs$ so $a\sim_{S_x}b$.
\end{xpl}

Another example to think of would be inverse semigroups, where $a\sim_Ab$ is just saying that $a$ and $b$ have a common lower bound in $A$.  We prove something similar for more general semigroups in \eqref{ZAAZ} below by considering the following subsets.
\begin{align*}
A^Z&=\{z\in Z:\exists a\in A\ (az=a)\}.\\
{}^Z\!A&=\{z\in Z:\exists a\in A\ (za=a))\}.
\end{align*}
First, however, we need to examine their basic properties.

\begin{prp}\label{A<A*}
If $A$ is an atlas then $A^Z$ is a subsemigroup of $Z$.  Moreover,
\begin{align}
\label{ZSources}A^Z=A^{<Z}={}^Z(A^*)&=(A^Z)^Z=\mathsf{s}(A)^Z\subseteq\mathsf{s}(A).\\
\label{Binormality}A\ni a''<_{a'}a\qquad&\Rightarrow\qquad a'{}^Z\!Aa\subseteq A^Z\\
\label{BoundedNormality}a\in A\cap S^>\qquad&\Rightarrow\qquad{}^Z\!Aa=aA^Z.
\end{align}
\end{prp}

\begin{proof}
If $y,z\in A^Z$ then we have $a,b\in A$ with $ay=a$ and $bz=b$.  As $A$ is an atlas, $A\subseteq A^<$ so we have $a',a''\in S$ with $A\ni a''<_{a'}a$.  Then $ba'a\in AA^*A\subseteq A$ and, moreover, $ba'ayz=ba'az=bza'a=ba'a$.  Thus $yz\in A^Z$, which shows that $A^Z$ is a subsemigroup of $Z$.  A dual argument works for ${}^Z\!A$.

\begin{itemize}
\item[\eqref{ZSources}]  If $z\in A^{*Z}$ then we have $a,a',a''\in S$ with $A\ni a''<_{a'}a$ and $a'z=a'$ and hence $za''=za''a'a=a''a'za=a''a'a=a''$, showing that $z\in{}^Z\!A$.  Conversely, if $z\in{}^Z\!A$, then we have $a\in A$ with $za=a$.  As $A\subseteq A^<$, we have $a',a''\in S$ with $A\ni a''<_{a'}a$ and hence $a'aa'z=a'zaa'=a'aa'\in A^*AA^*\subseteq A^*$, showing that $z\in A^{*Z}$.  Thus ${}^Z\!A=A^{*Z}$ and, likewise, ${}^Z(A^*)=A^Z$.  Then \eqref{**=} yields $A^{<Z}=A^{**Z}={}^Z(A^*)=A^Z$.

Now if $z\in A^Z$ then we have $a\in A$ with $az=a$.  Taking $a',a''\in A$ with $A\ni a''<_{a'}a$, we see that $A^Z\ni a'a=a'az$, so $z\in(A^Z)^Z$, and $a'a<_zz$, so $z\in(A^*A)^<=\mathsf{s}(A)$.  This shows that $A^Z\subseteq\mathsf{s}(A)$ and $A^Z\subseteq(A^Z)^Z\subseteq\mathsf{s}(A)^Z$.  Conversely, if $z\in\mathsf{s}(A)^Z$ then we have $b\in\mathsf{s}(A)$ with $bz=b$.  Taking $a\in A$, note $abz=ab\in A\mathsf{s}(A)\subseteq A^<$ so $z\in A^{<Z}=A^Z$, showing that $\mathsf{s}(A)^Z\subseteq A^Z$.

\item[\eqref{Binormality}] If $A\ni a''<_{a'}a$ and $zb=b\in A$ then we have $b',b''\in S$ with $A\ni b''<_{b'}b$ and hence $bb'a''\in AA^*A\subseteq A$ and $bb'a''a'za=zbb'a''a'a=bb'a''$, i.e. $bb'a''$ witnesses $a'za\in A^Z$.

\item[\eqref{BoundedNormality}] If $A\ni a<_{a'}a''$ and $z\in{}^Z\!A$ then \eqref{Binormality} yields $a'za''\in A^Z$ and hence $za=zaa'a''=aa'za''\in aA^Z$. \qedhere
\end{itemize}
\end{proof}

Dually to \eqref{Binormality}, $A\ni a''<_{a'}a$ also implies $aA^Za'\subseteq{}^Z\!A$.

\begin{prp}
If $A\subseteq S$ is an atlas and $a,b\in A^<$ then
\begin{align}
\label{ZAAZ}a\sim_Ab\qquad&\Leftrightarrow\qquad{}^Z\!AaA^Z\cap{}^Z\!AbA^Z\neq\emptyset\\
\label{ZAandAZ}&\Leftrightarrow\qquad\exists s\in A^*\ (sa=sb\in A^Z\text{ and }as=bs\in{}^Z\!A).
\end{align}
\end{prp}

\begin{proof}
Take $a,b\in A^<$.  If $s\in A^*$ and $sa=sb$ or $as=bs$ then $sas=sbs$ so $a\sim_Ab$.  This proves the last $\Leftarrow$.

If $a\sim_Ab$ then we have $s,t\in A^*$ with $sat=sbt$.  Taking any $s',s'',t',t''\in S$ with $A^*\ni s''<_{s'}s$ and $A^*\ni t''<_{t'}t$, we see that ${}^Z\!AaA^Z\ni s'satt'=s'sbtt'\in{}^Z\!AbA^Z$.  this proves the first $\Rightarrow$.

Conversely, assume ${}^Z\!AaA^Z\cap{}^Z\!AbA^Z\neq\emptyset$.  First we claim that $aA^Z\cap bA^Z\neq\emptyset$.  To see this, take $w,x\in{}^Z\!A$ and $y,z\in A^Z$ with $way=xbz$.  As $a\in A^<$, we have $a',a''\in S$ with $A\ni a''<_{a'}a$.  Note that $waya'a=waa'ay=aa'way\in aa'{}^Z\!AaA^Z\subseteq aA^Z$, by \autoref{A<A*}.  Likewise, as $b\in A^<$, we have $b',b''\in S$ with $A\ni b''<_{b'}b$ and then $xbzb'b\in bA^Z$.  As $way=xbz$,
\[aA^Z\supseteq aA^Zb'b\ni waya'ab'b=xbzb'ba'a\in bA^Za'a\subseteq bA^Z.\]
This proves the claim.

Next we claim that we have $t\in A^*$ with $at=bt\in{}^Z\!A$ and $ta\in A^Z$.  To see this, take $p,q\in A^Z$ with $ap=bq$.  As $p\in A^Z$, we have $c\in A\subseteq A^<$ with $cp=c$ and $c',c''\in S$ with $A\ni c''<_{c'}c$.  Likewise, we have $d,d',d''\in S$ with $A\ni d''<_{d'}d$ and $dq=d$.  Let $z=c'cd'd=d'dc'c\in A^Z\cap A^*A$, by \autoref{A<A*}.  Now just note that $az=ac'cd'd=ac'cpd'd=apc'cd'd=apz$ and, likewise,
\[bz=bd'dc'c=bd'dqc'c=bqd'dc'c=bqz=apz=az.\]
Taking the $a'$ above, we see that $t=za'\in A^*AA^*\subseteq A^*$, $ta=za'a\in A^ZA^Z\subseteq A^Z$ and $bza'=aza'\in aA^Za'\subseteq{}^Z\!A$, by \autoref{A<A*}, proving the claim.

Likewise, we have $u\in A^*$ with $ua=ub\in A^Z$ and $au\in{}^Z\!A$.  \autoref{A<A*} again yields $taub=taua\in A^ZA^Z\subseteq A^Z$ and $btau=atau\in A^ZA^Z\subseteq A^Z$, i.e.
\[sa=sb\in A^Z\text{ and }as=bs\in{}^Z\!A,\]
for $s=tau$.  This completes the cycle of equivalences.
\end{proof}

One could also write down other characterisations of $\sim_A$ that lie somewhere in between \eqref{ZAAZ} and \eqref{ZAandAZ}.  One could also replace $A^*$ with some coinitial subset $A'$.  That is to say if $A$ is an atlas and $A'\subseteq A^*=(A')^<$ then, for example,
\begin{equation}\label{Aprime}
a\sim_Ab\qquad\Leftrightarrow\qquad\exists s\in A'\ (as=bs).
\end{equation}
Indeed, if we have $s\in A'\subseteq A^*$ with $as=bs$ then, in particular, $aA^Z\cap bA^Z\neq\emptyset$ and so certainly ${}^Z\!AaA^Z\cap{}^Z\!AbA^Z\neq\emptyset$ and hence $a\sim_Ab$, by \eqref{ZAAZ}.  Conversely, if $a\sim_Ab$ then we have $s\in A^*$ with $as=bs$, by \eqref{ZAandAZ}.  If $A^*=(A')^<$ too then we have $a',t\in S$ with $A'\ni a'<_ts$ and hence $aa'=asta'=bsta'=ba'$, proving \eqref{Aprime}.

Now we can finally show that $\sim_A$ is an equivalence relation.

\begin{cor}
If $A\subseteq S$ is an atlas then $\sim_A$ is an equivalence relation on $A^<$.
\end{cor}

\begin{proof}
We immediately see that $\sim_A$ is a reflexive and symmetric relation on $A^<$.  On the other hand, if $a\sim_Ab\sim_Ac$ then we have $s,t\in A^*$ with $sa=sb$ and $bt=ct$, by \eqref{ZAandAZ}.  Thus $sat=sbt=sct$, showing that $a\sim_Ac$.  This shows that $\sim_A$ is also transitive and hence an equivalence relation on $A^<$.
\end{proof}

Also, multiples of ${}^Z\!A$ and $A^Z$ do not change the equivalence class.

\begin{cor}
For any atlas $A\subseteq S$, $y\in{}^Z\!A$, $z\in A^Z$ and $a\in A^<$,
\begin{equation}\label{azsimAa}
ya\sim_Aa\sim_Aaz.
\end{equation}
\end{cor}

\begin{proof}
As $y^2\in{}^Z\!A$, $z^2\in A^Z$ and $y(ya)z^2=y^2az^2=y^2(az)z$, \eqref{ZAAZ} yields \eqref{azsimAa}.
\end{proof}

Moreover, these equivalence relations even respect products.

\begin{cor}\label{abAB}
If $A,B\subseteq S$ are atlases with $\mathsf{s}(A)=\mathsf{r}(B)$ then
\begin{equation}\label{abAB}
a\sim_Ac\quad\text{and}\quad b\sim_Bd\qquad\Rightarrow\qquad ab\sim_{AB}cd.
\end{equation}
\end{cor}

\begin{proof}
If $a\sim_Ac$ and $b\sim_Bd$ then \eqref{ZAandAZ} yields $s\in A^*$ and $t\in B^*$ with $sa=sc$ and $bt=dt$.  Then $tsabts=tscdts$ so $ab\sim_{AB}cd$, as $ts\in B^*A^*\subseteq(AB)^*$ (noting $AB$ is itself an atlas because $\mathsf{s}(A)=\mathsf{r}(B)$ -- see \autoref{ABatlas}).
\end{proof}

\section{Bundles}\label{Bundles}

In the present section we show that the equivalence classes of the previous section form an \'etale bundle over the coset groupoid.

Firstly, if $A\subseteq S$ is an atlas and $a\in A^<$, we define $[a,A]$ to be the pair consisting of the $\sim_A$-equivalence class of $a$ together with coset $A^<$ generated by $A$, i.e.
\[[a,A]=(a^{\sim_A},A^<).\]
We denote the set of all such pairs by
\[\widetilde{\mathcal{C}}(S)=\{[c,C]:c\in C\in\mathcal{C}(S)\}.\]

\begin{prp}
$\widetilde{\mathcal{C}}(S)$ is a groupoid under the product
\[[a,A][b,B]=[ab,AB]\quad\text{when}\quad\mathsf{s}(A)=\mathsf{r}(B),\]
where $[c,C]^{-1}=[c',C^*]$, for any $c'$ with $C\ni d<_{c'}c$.
\end{prp}

\begin{proof}
By \eqref{abAB}, the product is well-defined.  Also, as $S$ is a semigroup and $\mathcal{C}(S)$ is a groupoid, $\widetilde{\mathcal{C}}(S)$ can be viewed as a semicategory with unit cosets as objects.  To show $\widetilde{C}(S)$ is a category we must find a unit $[z,U]\in\widetilde{\mathcal{C}}(S)$ for every unit $U\in\mathcal{C}(S)$.

We claim $[z,U]\in\widetilde{\mathcal{C}}(S)$ is a unit, for any $z\in U^Z=U^{*Z}={}^ZU$, by \eqref{ZSources}.  Indeed, if $C\in\mathcal{C}(S)$ and $\mathsf{s}(C)=U$ then $C^Z=U^Z$, by \eqref{ZSources}.  For any $c\in C$, this means $cz\sim_Cc$ so $[c,C][z,U]=[cz,CU]=[c,C]$.  Likewise, if $\mathsf{r}(C)=U$ then ${}^ZC={}^ZU$ so $[z,U][c,C]=[zc,UC]=[c,C]$, proving the claim.  Thus $\widetilde{\mathcal{C}}(S)$ is a category.

Next note that, for any $c\in C$, we have $c',d\in S$ with $C\ni d<_{c'}c$ and then $[c',C^*][c,C]=[c'c,\mathsf{s}(C)]$ and $[c,C][c',C^*]=[cc',\mathsf{r}(C)]$.  As $c'c\in C^Z=\mathsf{s}(C)^Z$ and $cc'\in{}^ZC={}^Z\mathsf{r}(C)$, these are both units in $\widetilde{\mathcal{C}}(S)$ by what we just proved and hence
\begin{equation}\label{sCinverse}
[c',C^*]=[c,C]^{-1}.
\end{equation}
This shows that every $[s,C]\in\widetilde{\mathcal{C}}(S)$ is invertible so $\widetilde{\mathcal{C}}(S)$ is a groupoid.
\end{proof}

Let $\rho_\mathcal{C}:\widetilde{\mathcal{C}}(S)\rightarrow\mathcal{C}(S)$ denote the canonical projection onto $\mathcal{C}(S)$, i.e.
\[\rho_\mathcal{C}([a,A])=A^<.\]
The canonical topology on $\widetilde{\mathcal{C}}(S)$ is generated by the sets
\begin{align*}
\widetilde{\mathcal{C}}^s&=\{[s,C]:s\in C\in\mathcal{C}(S)\}\quad\text{and}\\
\widetilde{\mathcal{C}}_s&=\{[c,C]:c,s\in C\in\mathcal{C}(S)\}=\rho_\mathcal{C}^{-1}[\mathcal{C}_s].
\end{align*}

\begin{thm}\label{CosetBundleEtale}
$\rho_\mathcal{C}:\widetilde{\mathcal{C}}(S)\rightarrow\mathcal{C}(S)$ is an \'etale bundle.
\end{thm}

\begin{proof}
To see that the inverse map on $\widetilde{\mathcal{C}}(S)$ is continuous, say $[a,A]^{-1}\in\widetilde{\mathcal{C}}^{b'}$.  This means $[a,A]^{-1}=[b',A^*]$ so, in particular, $b'\in A^*$.  Taking $b,c\in S$ with $A\ni c<_{b'}b$, \eqref{sCinverse} yields $[a,A]=[b',A^*]^{-1}=[b,A]$ so $[a,A]\in\widetilde{\mathcal{C}}^b\cap\widetilde{\mathcal{C}}_c$.  Moreover, for any other $[b,B]\in\widetilde{\mathcal{C}}^b\cap\widetilde{\mathcal{C}}_c$, we see that $B\ni c<_{b'}b$ and hence $[b,B]^{-1}=[b',B^*]\in\widetilde{\mathcal{C}}^{b'}$, by \eqref{sCinverse}.  Likewise, as in the proof of \autoref{EtaleCosets}, if $[a,A]^{-1}\in\widetilde{\mathcal{C}}_{b'}$ then again $b'\in A^*$ so, taking $b,c\in S$ with $A\ni c<_{b'}b$, we see that $[a,A]\in\widetilde{\mathcal{C}}_c$ and $(\widetilde{\mathcal{C}}_c)^{-1}\subseteq\widetilde{\mathcal{C}}_{b'}$.  This shows that the inverse map is continuous.

For the product, say $[a,A][b,B]\in\widetilde{\mathcal{C}}^s$.  As $[a,A][b,B]=[ab,AB]$, it follows that $ab\sim_{AB}s$, so \eqref{Aprime} yields $a'\in A^*$ and $b'\in B^*$ with $abb'a'=sb'a'$.  Taking $a'',b'',c,d\in S$ with $A\ni c<_{a'}a''$ and $B\ni d<_{b'}b''$, we see that $[a,A]\in\widetilde{\mathcal{C}}^a\cap\widetilde{\mathcal{C}}_c$ and $[b,B]\in\widetilde{\mathcal{C}}^b\cap\widetilde{\mathcal{C}}_d$.  Moreover, for any other $[a,C]\in\widetilde{\mathcal{C}}^a\cap\widetilde{\mathcal{C}}_c$ and $[b,D]\in\widetilde{\mathcal{C}}^b\cap\widetilde{\mathcal{C}}_d$ such that $[a,C][b,D]$ is defined, \eqref{Aprime} again yields $ab\sim_{CD}s$, as $c\in C$ and $d\in D$ so $a'\in C^*$ and $b'\in D^*$, and hence $[a,C][b,D]=[ab,CD]\in\widetilde{\mathcal{C}}^s$.  Similarly, if $[a,A][b,B]\in\widetilde{\mathcal{C}}_s$, then we have $c\in A$ and $d\in B$ with $cd<s$.  Then $[a,A]\in\widetilde{\mathcal{C}}_c$ and $[b,B]\in\widetilde{\mathcal{C}}_d$.  Moreover, for any other $[e,C]\in\widetilde{\mathcal{C}}_c$ and $[f,D]\in\widetilde{\mathcal{C}}_d$ such that $[e,C][f,D]$ is defined, $[e,C][f,D]=[ef,CD]\in\widetilde{\mathcal{C}}_{cd}\subseteq\widetilde{\mathcal{C}}_s$.  This shows that the product is continuous.

So $\widetilde{\mathcal{C}}(S)$ is a topological groupoid, and we already know that $\mathcal{C}(S)$ is an \'etale groupoid, by \autoref{EtaleCosets}.  It only remains to prove that $\rho_\mathcal{C}$ is a locally injective continuous open isocofibration.  Local injectivity is immediate from the fact that $\rho_\mathcal{C}$ is injective on $\widetilde{\mathcal{C}}^s$, for all $s\in S$.  Continuity is immediate from the fact that $\widetilde{\mathcal{C}}_s=\rho_\mathcal{C}^{-1}[\mathcal{C}_s]$ is open, for all $s\in S$.  The definition of the product on $\mathcal{C}(S)$ also immediately shows that $\rho_\mathcal{C}$ is an isocofibration.

To see that $\rho_\mathcal{C}$ is also open, say $F$ and $G$ and finite subsets of $S$ and
\[[a,A]\in\widetilde{\mathcal{C}}^F_G=\bigcap_{f\in F}\widetilde{\mathcal{C}}^f\cap\bigcap_{g\in G}\widetilde{\mathcal{C}}_g.\]
For all $f\in F$, this means $a\sim_Af$ so we have $a'_f\in A^*$ with $aa'_f=fa'_f$.  Then we can take $a_f,b_f\in S$ with $A\ni b_f<_{a'_f}a_f$.  Note that
\[H=\{a\}\cup\{b_f:f\in F\}\cup F\cup G\subseteq A\]
and hence $\rho_\mathcal{C}([a,A])=A\in\mathcal{C}_H$.  Moreover, for all $B\in\mathcal{C}_H$ and $f\in F$, note $b_f\in B$ so $a'_f\in B^*$ and hence $a\sim_Bf$.  Thus $[a,B]\in\widetilde{\mathcal{C}}^F_G$ so $B=\rho_\mathcal{C}([a,B])\in\rho_\mathcal{C}[\widetilde{\mathcal{C}}^F_G]$, showing that $\mathcal{C}_H\subseteq\rho_\mathcal{C}[\widetilde{\mathcal{C}}^F_G]$.  This shows that $\rho_\mathcal{C}$ is also an open map.
\end{proof}

We refer to $\rho_\mathcal{C}:\widetilde{\mathcal{C}}(S)\rightarrow\mathcal{C}(S)$ as the \emph{coset bundle}.

\section{Representations}\label{TheUniversalRepresentation}

Recall the notion of an \'etale representation from \autoref{EtaleRep}.

\begin{dfn}\label{BundleRep}
A \emph{bundle representation} of the structured semigroup $(S,N,Z)$ on a groupoid bundle $\pi:F\rightarrow G$ is a semigroup homomorphism $\theta:S\rightarrow\mathcal{S}(\pi)$ to the slice-sections of $\pi$ such that $a\mapsto\mathrm{dom}(\theta(a))$ is an \'etale representation.
\end{dfn}

If $\theta:S\rightarrow\mathcal{S}(\pi)$ is a semigroup homomorphism then so is $a\mapsto\mathrm{dom}(\theta(a))$.  To verify it is an \'etale representation, it thus suffices to check that these domains cover $G$, \eqref{RoundRep} holds and $N$ gets mapped to $\mathcal{O}(G^0)$, i.e. $\theta[N]\subseteq\mathcal{N}(\pi)$.  Note that \eqref{RoundRep} here ensures that $\theta[S]$ is a local-inverse semigroup, as per \autoref{local} -- see \autoref{<Equiv}.

Again, we call a bundle representation $\mu:S\rightarrow\mathcal{S}(\pi)$ \emph{universal} if, for every bundle representation $\theta:S\rightarrow\mathcal{S}(\pi')$, there exists a unique Pierce morphism $(\phi,\tau)$ from $\pi$ to $\pi'$ such that $\theta=\tfrac{\tau}{\phi}\circ\mu$ for $\tfrac{\tau}{\phi}:\mathcal{S}(\pi)\rightarrow\mathcal{S}(\pi')$ defined as in \autoref{BundleMorphism->SemigroupHomo}.

\begin{thm}\label{CosetBundleUniversal}
If $Z$ is symmetric then $a\mapsto\widetilde{a}$ is a universal bundle representation where $\widetilde{a}\in\mathcal{S}(\rho_\mathcal{C})$ denotes the slice-section such that, for all $C\in\mathrm{dom}(\widetilde{a})=\mathcal{C}_a$,
\[\widetilde{a}(C)=[a,C].\]
\end{thm}

\begin{proof}
For any $a\in S$, $\mathrm{dom}(\widetilde{a})=\mathcal{C}_a$ is a slice, by \autoref{CaSlice}.  To see that $\widetilde{a}$ is continuous, first note $\widetilde{a}^{-1}[\widetilde{\mathcal{C}}_s]=\mathcal{C}_a\cap\mathcal{C}_s$.  On the other hand, if $\widetilde{a}(C)=[a,C]\in\widetilde{\mathcal{C}}^s$ then $az=sz$, for some $z\in C^Z$.  Take any $b\in C$ with $bz=b$.  If $B\in\mathcal{C}_a\cap\mathcal{C}_b$ then $z\in B^Z$ and hence $\widetilde{a}(B)=[a,B]=[s,B]\in\widetilde{\mathcal{C}}^s$.  Thus $C\in\mathcal{C}_a\cap\mathcal{C}_b\subseteq\widetilde{a}^{-1}[\widetilde{\mathcal{C}}^s]$.  We have shown that preimages of subbasic sets are open and hence $\widetilde{a}$ is continuous.

If $Z$ is symmetric then $a\mapsto\mathrm{dom}(\widetilde{a})=\mathcal{C}_a$ is an \'etale representation, thanks to \autoref{SymmetricEtaleRepresentation}.  In particular, for any $a,b\in S$, $\mathrm{dom}(\widetilde{\,ab\,})=\mathrm{dom}(\widetilde{a})\cdot\mathrm{dom}(\widetilde{b})$ and, for any $A\in\mathrm{dom}(\widetilde{a})$ and $B\in\mathrm{dom}(\widetilde{b})$ with $\mathsf{s}(A)=\mathsf{r}(B)$,
\[\widetilde{\,ab\,}(A\cdot B)=[ab,AB]=[a,A][b,B]=\widetilde{a}(A)\widetilde{b}(B).\]
This shows that $\widetilde{\,ab\,}=\widetilde{a}\widetilde{b}$, which in turn shows that $a\mapsto\widetilde{a}$ is a semigroup homomorphism and hence a bundle representation.

To see that $a\mapsto\widetilde{a}$ is universal, take another bundle representation $\theta:S\rightarrow\mathcal{S}(\pi)$.  By \autoref{SymmetricEtaleRepresentation}, we have a unique \'etale morphism $\phi:G\rightarrow\mathcal{C}(S)$ such that $\mathrm{dom}(\theta(a))=\phi^{-1}[\mathcal{C}_a]$, for all $a\in S$.  Indeed, the proof shows that, for all $g\in G$,
\[\phi(g)=\{a\in S:g\in\mathrm{dom}(\theta(a))\}.\]

We claim we have a continuous functor $\tau:\phi^{\rho_\mathcal{C}}\widetilde{\mathcal{C}}(S)\rightarrow F$ given by
\[\tau([a,\phi(g)],g)=\theta(a)(g).\]
To see that $\tau$ is well-defined, take $a,b\in\phi(g)$ with $a\sim_{\phi(g)}b$.  Then $as=bs$, for some $s\in \phi(g)^*=\phi(g^{-1})$, and hence
\[\theta(a)(g)\theta(s)(g^{-1})=(\theta(a)\theta(s))(\mathsf{r}(g))=\theta(as)(\mathsf{r}(g))=\theta(bs)(\mathsf{r}(g))=\theta(b)(g)\theta(s)(g^{-1}).\]
Multiplying by $\theta(s)(g^{-1})^{-1}$ yields $\theta(a)(g)=\theta(b)(g)$, as required.  To see that $\tau$ is a functor, note that  $([a,\phi(g)],g),([b,\phi(h)],h)\in\phi^{\rho_\mathcal{C}}\widetilde{\mathcal{C}}(S)$ and $\mathsf{s}(g)=\mathsf{r}(h)$ implies
\begin{align*}
\tau(([a,\phi(g)],g)([b,\phi(h)],h))&=\tau([ab,\phi(gh)],gh)=\theta(ab)(gh)=\theta(a)(g)\theta(b)(h)\\
&=\tau([a,\phi(g)],g)\tau([b,\phi(h)],h).
\end{align*}
For continuity, note that if $\tau([a,\phi(g)],g)=\theta(a)(g)\in O$, for some open $O\subseteq F$, then the continuity of $\theta(a)$ yields open $O'\subseteq G$ with $g\in O'$ and $\theta(a)[O']\subseteq O$.  Then $([a,\phi(g)],g)\in\widetilde{C}^a\times O'$ and, for any other $([a,\phi(h)],h)\in(\widetilde{C}^a\times O')$, we see that
\[\tau([a,\phi(h)],h)=\theta(a)(h)\in\theta(a)[O']\subseteq O.\]
This shows that $\tau$ is continuous, completing the proof of the claim.  Further noting that $\pi(\tau([a,\phi(g)],g))=\pi(\theta(a)(g))=g=\rho_{\mathcal{C}\phi}([a,\phi(g)],g)$, i.e. $\pi\circ\tau=\rho_{\mathcal{C}\phi}$, we see that $(\phi,\tau)$ is a Pierce morphism from $\rho_\mathcal{C}$ to $\pi$.

For all $g\in G$ and $a\in\phi(g)$, we see that
\[\theta(a)(g)=\tau([a,\phi(g)],g)=\tau(\widetilde{a}(\phi(g)),g)=\tau((\widetilde{a}\circ\phi)(g),g)=\tfrac{\tau}{\phi}(\widetilde{a})(g),\]
where $\tfrac{\tau}{\phi}:\mathcal{S}(\rho_\mathcal{C})\rightarrow\mathcal{S}(\pi)$ is the homomorphism from \autoref{BundleMorphism->SemigroupHomo}.  Also
\[\mathrm{dom}(\tfrac{\tau}{\phi}(\widetilde{a}))=\phi^{-1}[\mathrm{dom}(\widetilde{a})]=\phi^{-1}[\mathcal{C}_a]=\mathrm{dom}(\theta(a))\]
and hence $\theta(a)=\tfrac{\tau}{\phi}(\widetilde{a})$, for all $a\in S$, i.e. $\theta=\tfrac{\tau}{\phi}\circ\widetilde{\mathcal{C}}$, where $\widetilde{\mathcal{C}}(a)=\widetilde{a}$.

To see that $(\phi,\tau)$ is the unique, take another Pierce morphism $(\phi',\tau')$ with $\theta=\tfrac{\tau'}{\phi'}\circ\widetilde{\mathcal{C}}$.  It follows that $a\mapsto\mathrm{dom}(\theta(a))=\mathrm{dom}(\tfrac{\tau'}{\phi'}(\widetilde{a}))=\phi'^{-1}[\mathcal{C}_a]$ is an \'etale morphism and hence $\phi'=\phi$, by the uniqueness part of \autoref{SymmetricEtaleRepresentation}.  Now, for any $([a,\phi(g)],g)\in\phi^{\rho_\mathcal{C}}\widetilde{\mathcal{C}}(S)$, just note that
\[\tau'([a,\phi(g)],g)=\tfrac{\tau'}{\phi'}(\widetilde{a})(g)=\theta(a)(g)=\tau([a,\phi(g)],g).\]
Thus $\tau'=\tau$, proving uniqueness and hence universality.
\end{proof}

It follows that cosets and their equivalence classes determine when structured semigroups have faithful bundle representations, at least under symmetry.

\begin{cor}\label{CosetBundleFaithful}
If $Z$ is symmetric then the following are equivalent.
\begin{enumerate}
\item\label{AnyFaithful} $(S,N,Z)$ has a faithful bundle representation.
\item\label{CosetFaithful} The coset bundle representation of $(S,N,Z)$ is faithful.
\item\label{CosetSim} For all distinct $a,b\in S$, we have $C\in\mathcal{C}(S)$ with $\{a,b\}\cap C\neq\emptyset$ and $a\not\sim_Cb$.
\end{enumerate}
\end{cor}

\begin{proof}
Certainly \eqref{CosetFaithful} implies \eqref{AnyFaithful}.  Conversely, if $\theta$ is a faithful bundle representation then \autoref{CosetBundleUniversal} yields a Pierce morphism $(\phi,\tau)$ with $\theta=\tfrac{\tau}{\phi}\circ\widetilde{\mathcal{C}}$, where $\widetilde{\mathcal{C}}$ is the coset bundle representation.  As $\theta$ is injective, so is $\widetilde{\mathcal{C}}$, showing that \eqref{AnyFaithful} implies \eqref{CosetFaithful}.

More explicitly, \eqref{CosetFaithful} is saying that, for any distinct $a,b\in S$, we have $C\in\mathcal{C}(S)$ such that either only one of $\widetilde{a}$ and $\widetilde{b}$ are defined at $C$, and so certainly $a\not\sim_Cb$, or they are both defined at $C$ but have different values at $C$, i.e.
\[[a,C]=\widetilde{a}(C)\neq\widetilde{b}(C)=[b,C],\]
which again means $a\not\sim_Cb$.  This is immediately seen to be equivalent to \eqref{CosetSim}.
\end{proof}

Unlike the situation with inverse semigroups in \autoref{InverseSemigroupFaithful}, it is indeed possible for a (symmetric) structured semigroup to have no faithful representations.

\begin{xpl}\label{NonFaithfulExample}
Take $S=N=Z=\{0,a\}$ where $a^2=a0=0a=00=0$.  This is certainly a structured semigroup with symmetric $Z$.  For any representation $a\mapsto\widehat{a}$ on a groupoid bundle $\pi:F\rightarrow G$, we would have $\mathrm{dom}(\widehat{a})=\mathrm{dom}(\widehat{0})\subseteq G^0$, as $a^2=0\in N$, and $\mathrm{ran}(\widehat{0})\subseteq F^0$, as $00=0$, and hence $\widehat{0}=\widehat{a0}=\widehat{a}\widehat{0}=\widehat{a}$, showing that the representation can not be faithful.
\end{xpl}

\begin{rmk}
On the one hand, this might motivate one to look for even more general kinds of representations than those arising from slice-sections of groupoid bundles.  On the other hand, one might look for some natural restrictions or additional structures to place on structured semigroups to ensure that they do have faithful bundle representations.  In our subsequent work (see \cite{Bice2021Pierce} and \cite{Bice2021DHKR}), we actually do a bit both, considering semigroups that sit nicely within larger rings and their representations as sections of more general category bundles.
\end{rmk}

\section{Filters}\label{Filters}

As usual, we call $D\subseteq S$ \emph{directed} if
\[\tag{Directed}a,b\in D\qquad\Rightarrow\qquad\exists d\in D\ (d<a,b).\]
A \emph{filter} is a directed up-set $D=D^<$.  In particular, the non-empty directed cosets
\[\mathcal{D}(S)=\{D\in\mathcal{C}(S):D\text{ is directed}\}\]
are all filters.  Directed cosets offer a number of advantages over arbitrary cosets.  For example, one notes immediately notes that the canonical subbasis of $\mathcal{C}(S)$ becomes a basis when restricted to $\mathcal{D}(S)$, namely $(\mathcal{D}_a)_{a\in S}$, where $\mathcal{D}_a=\mathcal{C}_a\cap\mathcal{D}(S)$.  Our goal here is to show that the directed coset subbundle is still a groupoid bundle and that the corresponding subrepresentation of $S$ is `Zakrzewski-universal'.

\begin{rmk}
There will be plenty of directed cosets arising from principal filters in our primary motivating examples, as already noted in \autoref{PrincipalFilters}.  But in general, a structured semigroup $S$ may have few directed cosets -- again we will see in \autoref{DirectedCosetBundleFaithful} below that we have faithful bundle representations precisely when there are enough directed cosets to distinguish the elements of $S$ via their corresponding equivalence relations, at least when $Z$ is symmetric.
\end{rmk}

First we look at products.  Recall the notation $c|D$ from \autoref{ActionNotation}.

\begin{prp}
For any $c\in S$ and $D\subseteq S$,
\begin{equation}\label{bCdirected}
c|D\text{ and }D\text{ is directed}\qquad\Rightarrow\qquad cD\text{ is directed}.
\end{equation}
\end{prp}

\begin{proof}
If $c'cd=d\in D$ and $D$ is directed then, for any $a,b\in D$, we have $e\in D$ with $e<a,b,d$.  In particular, we have $d'\in S$ with $e<_{d'}d$ and hence $e=dd'e=c'cdd'e=c'ce$.  Then \eqref{ZSplitting} (applied to $S^\mathrm{op}$) yields $ce<ca,cb$, showing that $cD$ is directed.
\end{proof}

Likewise, duals respect directedness.

\begin{prp}\label{*directed}
If $D\subseteq S$ is directed then so is $D^*$ and
\begin{equation}\label{D*DD*}
D^*DD^*\subseteq D^{*<}.
\end{equation}
\end{prp}

\begin{proof}
Take $b',d'\in D^*$, so we have $a,b,c,d\in S$ with $D\ni a<_{b'}b$ and $D\ni c<_{d'}d$.  As $D$ is directed, we have $e,f\in D$ with $e<_{f'}f< a,c$ so \eqref{Transitivity} yields $e,f<_{b'}b$ and $e,f<_{d'}d$.  Then \eqref{Switch} yields $b'fd'<_bb'$ and $b'fd'<_dd'$.  Also $e<_{b'fd'}df'b$, as $eb'fd'\in NN\subseteq N$ and $eb'fd'df'b=eb'ff'b=ff'eb'b=e$.  Thus $b'fd'\in D^*$, showing that $D^*$ is indeed directed.

Now say $e\in D$ and $b',d'\in D^*$, so we again have $a,b,c,d\in S$ with $D\ni a<_{b'}b$ and $D\ni c<_{d'}d$.  As $D$ is directed, we have $f,g\in D$ with $g<_{f'}f<a,c,e$.  We immediately verify $g<_{b'fd'}df'b$ (e.g. $gb'fd'df'b=gd'db'ff'b=ff'gb'b=g$) so $b'fd'\in D^*$.  Taking $e'\in S$ with $f<_{e'}e$, we also immediately verify $b'fd'<_{de'b}b'ed'$ (e.g. $b'fd'de'bb'ed'=b'fe'bb'ed'=b'bb'fe'ed'=b'fd'$) and hence $b'ed'\in D^{*<}$, showing that $D^*DD^*\subseteq D^{*<}$.
\end{proof}

It follows that, when $N$ is diagonal, the directed cosets are precisely the filters and, moreover, for any $(d^{\sim_D},D)\in\widetilde{\mathcal{D}}(S)$, the second coordinate is just the up-closure of the first (where $\widetilde{\mathcal{D}}(S)=\rho_\mathcal{C}^{-1}[\mathcal{D}(S)]=\{(d^{\sim_D},D)\in\widetilde{\mathcal{C}}(S):D\in\mathcal{D}(S)\}$).

\begin{prp}\label{FilterCosets}
Any directed $D=D^{**}\subseteq S$ is both a filter and a coset.\\
Conversely, if $N$ is diagonal then every filter is a coset and
\[d\in D\in\mathcal{D}(S)\qquad\Rightarrow\qquad D=(d^{\sim_D})^<.\]
\end{prp}

\begin{proof}
As $D$ is directed and hence round, $D^<\subseteq D^{<<}\subseteq D^{**}=D$, i.e. $D$ is also an up-set and hence a filter.  Also \eqref{Antimorphism}, \eqref{<*} and \eqref{D*DD*} yield
\[DD^*D=D^{**}D^*D^{**}\subseteq(D^*DD^*)^*\subseteq D^{*<*}\subseteq D^{**<}=D^<=D,\]
so $D$ is also a coset.

When $N$ is diagonal, any filter $D\subseteq S$ automatically satisfies $D=D^{<<}=D^{**}$, by \eqref{DiagonalStars}, and is thus a coset by the above argument.  Moreover, for any $c,d\in D$, we have $a,b',b\in S$ with $D\ni a<_{b'}b<c,d$.  Taking $c',d'\in D$ with $b<_{c'}c$ and $b<_{d'}d$, we see that $cc'bb'd=bb'd$ and $c'bb'd=c'dd'bb'd\in N$, by diagonality, because $c'dd'b,d'b,d'bb'd\in N$.  Thus $d\sim_Dbb'd<_{c'}c$, by \eqref{azsimAa}, i.e. $c\in(d^{\sim_D})^<$, showing that $D\subseteq(d^{\sim_D})^<\subseteq D^<\subseteq D$.
\end{proof}

\begin{dfn}
In a groupoid $G$, we call $I\subseteq G$ an \emph{ideal} if, for all $(g,h)\in G^2$,
\[\tag{Ideal}\label{Ideal}g\in I\quad\text{or}\quad h\in I\qquad\Rightarrow\qquad gh\in I.\]
\end{dfn}

In particular, if $I\subseteq G$ is an ideal then $g\in I$ implies $gg^{-1}\in I$ and hence $g^{-1}=g^{-1}gg^{-1}\in I$, i.e. $I=I^{-1}=II$.  It follows that ideals are precisely the replete/isomorphism-closed full subgroupoids.  Moreover, if $G$ has an \'etale topology then the subspace topology on any ideal is again \'etale \textendash\, see \cite{Bice2021}.

\begin{prp}\label{DIdeal}
$\mathcal{D}(S)$ is an ideal and hence an \'etale subgroupoid of $\mathcal{C}(S)$.
\end{prp}

\begin{proof}
If $C\in\mathcal{C}(S)$, $D\in\mathcal{D}(S)$ and $(C,D)\in\mathcal{C}^2$ then $cD$ is directed, by \eqref{bCdirected}, and hence so is $(CD)^*=(cD)^*$ and $(CD)^{**}$, by \autoref{*directed} and \eqref{BC=bC}.  Likewise, $(D,C)\in\mathcal{C}^2$ implies $(DC)^{**}$ is directed and hence $\mathcal{D}(S)$ an ideal.
\end{proof}

It will be useful to note units in $\mathcal{D}(S)$ are just units in $\mathcal{C}(S)$ generated by $N$.

\begin{prp}\label{DirectedUnitCosets}
Let $U\in\mathcal{C}(S)$ be a unit coset.
\begin{enumerate}
\item $(U\cap N)^<$ is a directed unit coset.
\item $(U\cap N)^<=U$ when $U$ is directed.
\item\label{DUCdiagonal} $(U\cap N)^<=U^{Z<}$ when $N$ is diagonal.
\end{enumerate}
\end{prp}

\begin{proof}\
\begin{enumerate}
\item First we claim $U\cap N$ is directed.  Indeed, if $m,n\in U\cap N$ then we have $U\ni t<_{m'}m$ and $U\ni u<_{n'}n$.  Noting that $t=(tm')m\in NN\subseteq N$ and $u=(un')n\in NN\subseteq N$, it follows that $tu\in NN\subseteq N$.  As $U$ is a unit coset, $tu\in UU=UU^*\subseteq U$.  Thus $U\cap N\ni tu<m,n$, by \eqref{NInvariance}.

This proves the claim, and hence $(U\cap N)^<$ and $(U\cap N)^{<**}$ are also directed, by \autoref{*directed}.  Moreover, by \autoref{Cunit},
\[\emptyset\neq U\cap N\subseteq(U\cap N)^<\subseteq(U\cap N)^{<<<}\subseteq(U\cap N)^{<**}.\]
Also, for any $u\in(U\cap N)^{<**}\subseteq U^{<**}=U$ and $n\in U\cap N\subseteq(U\cap N)^{<**}$, directedness yields $t\in(U\cap N)^{<**}\subseteq U$ with $t<u,n$.  Then $t\in N$ too, as $n\in N$, i.e. $t\in U\cap N$ and hence $u\in(U\cap N)^<$.  This proves that $(U\cap N)^<=(U\cap N)^{<**}$ and hence $(U\cap N)^<$ is a unit coset, by \autoref{Cunit} and \autoref{FilterCosets}.

\item In particular, if $U=(U\cap N)^<$ then $U$ is directed.  Conversely, if $U$ is directed then, as above, for any $u\in U$ and $n\in U\cap N$, we have $t\in U$ with $t<u,n$ and hence $t\in N$, as $n\in N$.  Thus $u\in(U\cap N)^<$, showing that $U\subseteq(U\cap N)^<\subseteq U^<=U$.

\item By \eqref{ZSources}, $U^Z\subseteq\mathsf{s}(U)\cap Z\subseteq U\cap N$ and hence $U^{Z<}\subseteq(U\cap N)^<$.  For the reverse inclusion, assume $N$ is diagonal and take $u\in(U\cap N)^<$ so we have $t,u'\in S$ with $U\cap N\ni t<_{u'}u$.  Then we can further take $s,t'\in S$ with $U\ni s<_{t'}t$.  As $t't,t,tu'\in N$ and $N$ is diagonal, $t'tu'\in N$ and hence $U^Z\ni t't<_{u'}u$, showing that $(U\cap N)^<\subseteq U^{Z<}$. \qedhere
\end{enumerate}
\end{proof}

While $\mathcal{D}(S)$ might have a nice basis, it still need not be Hausdorff or even $T_1$ (when $C\subsetneqq D$, for $C,D\in\mathcal{D}(S)$, any neighbourhood of $C$ is a neighbourhood of $D$).  For better separation properties, we can consider ultrafilters
\[\mathcal{U}(S)=\{U\subseteq S:U\text{ is a maximal proper filter}\},\]
at least when $Z$ has a zero, i.e. $0s=s0=0\in Z$, for all $s\in S$.

\begin{thm}\label{IdealUltrafilters}
If $Z$ has a zero then $\mathcal{U}(S)$ is a locally Hausdorff ideal of $\mathcal{C}(S)$.
\end{thm}

\begin{proof}
As $0\in Z\subseteq N$, $0<_0s$, for all $s\in S$.  So a filter $F\subseteq S$ is proper iff $0\notin F$.

First we show that every $U\in\mathcal{U}(S)$ is a coset.  By \eqref{**}, $U=U^{<<}\subseteq U^{**}$.  As $0\notin U$, it follows that $0\notin U^*$, $0\notin U^{**}$ and hence $0\notin U^{**<}$.  By \autoref{*directed}, $U^*$ and $U^{**}$ are directed and hence $U^{**<}$ is a proper filter containing $U^{**}$ and hence $U$.  As $U$ is an ultrafilter, $U=U^{**}=U^{**<}$ and hence $U$ is a coset, by \autoref{FilterCosets}.

Now say $C$ is another coset with $(C,U)\in\mathcal{C}^2$ and take any $c,d\in C$ and $c'\in C^*$ with $cc'd=d$.  We already saw above that $(CU)^{**}=(cU)^{**}$ is a filter.  Extending to an ultrafilter $V\supseteq(cU)^{**}$, we then see that $(c'V)^{**}$ is another filter containing $U$ and hence $U=(c'V)^{**}$ by maximality.  Then $(cU)^{**}=(cc'V)^{**}=V$, i.e. $(CU)^{**}=(cU)^{**}$ was already maximal.  Likewise $(U,C)\in\mathcal{C}^2$ implies $(UC)^{**}$ is an ultrafilter, showing that the ultrafilters indeed form an ideal in $\mathcal{C}(S)$.

It follows that $\mathcal{U}(S)$ is an \'etale subgroupoid of $\mathcal{C}(S)$, by \cite[Proposition 2.7]{Bice2021}.  In particular, $\mathcal{U}(S)$ is locally homeomorphic to its unit space so to prove that $\mathcal{U}(S)$ is locally Hausdorff, it suffices to show that its unit space is Hausdorff.  To see this, take any distinct unit ultrafilters $U,V\in\mathcal{U}(S)$.  By \eqref{Multiplicativity}, $(UV)^<$ is again a filter.  By \autoref{Cunit}, we have $m\in U\cap N$ and $n\in V\cap N$.  By \eqref{NInvariance}, $U\subseteq(Un)^<\subseteq(UV)^<$ and $V\subseteq(mV)^<\subseteq(UV)^<$.  As $U$ and $V$ are distinct and maximal, we must have $(UV)^<=S$.  In particular, we must have $u\in U$ and $v\in V$ with $uv<0$ and hence $uv=0$.  Thus $\mathcal{U}_u\cap\mathcal{U}_v=\mathcal{C}_u\cap\mathcal{C}_v\cap\mathcal{U}(S)=\emptyset$.  Indeed, if we had $W\in\mathcal{U}_u\cap\mathcal{U}_v$ then $0=uv\in WW=WW^*\subseteq\mathsf{r}(W)=W$, as $W$ has to be a unit coset, again by \autoref{Cunit}.  But this would mean $W=0^<=S$ and, in particular, $W\notin\mathcal{U}(S)$.  This shows that the unit space of $\mathcal{U}(S)$ is Hausdorff and hence the entirety of $\mathcal{U}(S)$ is locally Hausdorff.
\end{proof}

\begin{rmk}
We will have more to say about $\mathcal{U}(S)$ in future work.  For the moment we simply remark that ultrafilters are more natural to consider than arbitrary filters when dealing with groupoid C*-algebras and more general `bumpy semigroups' \textendash\, see \cite{Bice2021} and \cite{BiceClark2021}.
\end{rmk}

\subsection{Zakrzewski-Universality}\label{Universality}

By \autoref{DIdeal}, the identity embedding of $\mathcal{D}(S)$ in $\mathcal{C}(S)$ is an \'etale functor.  For Zakrzewski-universality, we need a Zakrzewski morphism in the opposite direction.

\begin{dfn}
For any $C\in\mathcal{C}(S)$, we define
\[D\vartriangleleft C\qquad\Leftrightarrow\qquad D\text{ is a maximal directed subset of }C.\]
\end{dfn}
So $\mathrm{dom}(\vartriangleleft)\subseteq\mathcal{C}(S)$ and $\mathrm{ran}(\vartriangleleft)\subseteq\mathcal{P}(S)$.  In fact, $\mathrm{ran}(\vartriangleleft)\subseteq\mathcal{D}(S)$.

\begin{prp}\label{tri}
If $D\vartriangleleft C$ then $D\in\mathcal{D}(S)$.  If $\mathcal{D}(S)\ni D\subseteq C\in\mathcal{C}(S)$ then
\begin{equation}\label{triEq}
\mathsf{r}(C)\cap N=\mathsf{r}(D)\cap N\,\quad\,\Leftrightarrow\,\quad\,D\vartriangleleft C\,\quad\,\Leftrightarrow\,\quad\,\mathsf{s}(C)\cap N=\mathsf{s}(D)\cap N.
\end{equation}
Moreover, if $d\in C\in\mathcal{C}(S)$ then we have a unique $D\in\mathcal{D}(S)$ with $d\in D\vartriangleleft C$.
\end{prp}

\begin{proof}
If $D\vartriangleleft C$ then $C\in\mathcal{C}(S)$ and, in particular, $C=C^{**}$.  As $D$ is directed, so too are $D^*$ and $D^{**}$, by \autoref{*directed}.  Then \eqref{**} yields
\[D\subseteq D^<\subseteq D^{<<}\subseteq D^{**}\subseteq C^{**}\subseteq C\]
so $D=D^{**}$, by maximality.  Thus $D$ is a coset, by \autoref{FilterCosets}, i.e. $D\in\mathcal{D}(S)$.

Assume $d\in C\in\mathcal{C}(S)$.  Taking $c,d'\in C$ with $c<_{d'}d$, we see that $c=cd'd$ and $cd'\in\mathsf{r}(C)\cap N\subseteq(\mathsf{r}(C)\cap N)^<$ (as $\mathsf{r}(C)\cap N$ is directed \textendash\, see the proof of \autoref{DirectedUnitCosets}).  Thus $(\mathsf{r}(C)\cap N)^<|d$ so, by \autoref{A|b}, we have a coset
\[D=((\mathsf{r}(C)\cap N)^<d)^<\subseteq(\mathsf{r}(C)C)^<=C\]
with $\mathsf{r}(D)=\mathsf{r}((\mathsf{r}(C)\cap N)^<)=(\mathsf{r}(C)\cap N)^<$, by \autoref{DirectedUnitCosets}.  Furthermore,
\begin{equation}\label{rCN}
\mathsf{r}(C)\cap N=\mathsf{r}(D)\cap N
\end{equation}
because $\mathsf{r}(C)\cap N\subseteq (\mathsf{r}(C)\cap N)^<\cap N\subseteq\mathsf{r}(C)\cap N$.  Also $(\mathsf{r}(C)\cap N)^<d\ni cd'd<_{d'}d$ and hence $d\in D\in\mathcal{D}(S)$, by \autoref{DirectedUnitCosets} and \eqref{bCdirected}.

We claim that $D$ is actually the largest element of $\mathcal{D}(S)$ with $d\in D\subseteq C$.  To see this, take any other $A\in\mathcal{D}(S)$ with $d\in A\subseteq C$.  Then $\mathsf{r}(A)$ is directed, by \autoref{DIdeal}, and hence $\mathsf{r}(A)=(\mathsf{r}(A)\cap N)^<\subseteq(\mathsf{r}(C)\cap N)^<=\mathsf{r}(D)$, by \autoref{DirectedUnitCosets}.  Then \eqref{BC=bC} yields
\[A=(\mathsf{r}(A)A)^<=(\mathsf{r}(A)d)^<\subseteq(\mathsf{r}(D)d)^<=(\mathsf{r}(D)D)^<=D,\]
proving the claim.  Moreover, note that if we had $\mathsf{r}(A)\cap N=\mathsf{r}(C)\cap N$ then we would have $\mathsf{r}(A)=(\mathsf{r}(A)\cap N)^<=(\mathsf{r}(C)\cap N)^<=\mathsf{r}(D)$ which would yield $A=D$ above.  It follows that $D$ is the unique element of $\mathcal{D}(S)$ with $d\in D\vartriangleleft C$ and also the unique element of $\mathcal{D}(S)$ with $d\in D\subseteq C$ and $\mathsf{r}(C)\cap N=\mathsf{r}(D)\cap N$, which yields the first $\Leftrightarrow$ in \eqref{triEq}.  The second follows by a dual argument.
\end{proof}

Note that the second statement above is saying that any coset $C$ can be partitioned uniquely into maximal directed sub(co)sets.  In fact, uniqueness holds even among more general families of directed cosets covering $C$.

\begin{prp}\label{UniqueCover}
For any $C\in\mathcal{C}(S)$,
\[C^\vartriangleright=\{D\in\mathcal{D}(S):D\vartriangleleft C\}\]
is the unique subfamily of $\mathcal{D}(S)$ with union $C$ and a common source (or range).
\end{prp}

\begin{proof}
Take a family $\mathscr{D}\subseteq\mathcal{D}(S)$ with $C=\bigcup\mathscr{D}$.  We claim that
\begin{equation}\label{sourceUnion}
\mathsf{s}(C)\cap N=\bigcup_{D\in\mathscr{D}}\mathsf{s}(D)\cap N.
\end{equation}
To see this, take any $n\in\mathsf{s}(C)\cap N$.  By \autoref{tri}, we have $B\vartriangleleft C$ with $\mathsf{s}(B)\cap N=\mathsf{s}(C)\cap N$ so $n\in\mathsf{s}(B)=(B^*B)^<$, i.e. we have $b\in B$ and $b'\in B^*$ such that $b'b<n$.  Taking $e,f\in B$ with $e<_{b'}f$, directedness yields $a\in B$ with $a<b,e$ and hence $a<_{b'}f$.  As $C=\bigcup\mathscr{D}$, we have $D\in\mathscr{D}$ containing $a$ so $b\in D$ and $b'\in D^*$, which means $n\in(D^*D)^<\cap N=\mathsf{s}(D)\cap N$, proving the claim.

If the directed cosets in $\mathscr{D}$ also have a common source then the only way \eqref{sourceUnion} could hold is if $\mathsf{s}(C)\cap N=\mathsf{s}(D)\cap N$, for all $D\in\mathscr{D}$.  Then \autoref{tri} again yields $\mathscr{D}=C^\vartriangleright$, as required.
\end{proof}

When $N$ is diagonal, we can also determine when $D\vartriangleleft C$ from $C^Z$ and $D^Z$.

\begin{prp}
For any $C\in\mathcal{C}(S)$ and $D\in\mathcal{D}(S)$,
\begin{equation}\label{CtriD=>CZ=DZ}
D\vartriangleleft C\qquad\Rightarrow\qquad C^Z=D^Z.
\end{equation}
The converse also holds when $D\subseteq C$ and $N$ is diagonal.
\end{prp}

\begin{proof}
First note that $C^Z=(\mathsf{s}(C)\cap N)^Z$ because, by \eqref{ZSources},
\[C^Z=(C^Z)^Z\subseteq(\mathsf{s}(C)\cap N)^Z\subseteq\mathsf{s}(C)^Z=C^Z.\]
Thus, by \eqref{triEq}, $D\vartriangleleft C$ implies $\mathsf{s}(C)\cap N=\mathsf{s}(D)\cap N$ and hence $C^Z=D^Z$.

Conversely, if $N$ is diagonal, $D\subseteq C$ and $C^Z=D^Z$ then \autoref{DirectedUnitCosets} yields
\[\mathsf{s}(C)\cap N\subseteq(\mathsf{s}(C)\cap N)^<=C^{Z<}=D^{Z<}=(\mathsf{s}(D)\cap N)^<=\mathsf{s}(D).\]
Thus $\mathsf{s}(C)\cap N=\mathsf{s}(D)\cap N$ and hence $D\vartriangleleft C$, by \eqref{triEq}.
\end{proof}

Next we show that $\vartriangleleft$ is a Zakrzewski morphism such that the coset representation $\mathcal{C}$ factors through the directed coset representation $\mathcal{D}$ via preimages of $\vartriangleleft$.

\begin{thm}\label{triEtaleMorphism}
$\vartriangleleft$ is a Zakrzewski morphism from $\mathcal{C}(S)$ to $\mathcal{D}(S)$ with
\begin{equation}\label{CtriD}
\mathcal{C}=\overline{\vartriangleleft}\circ\mathcal{D}.
\end{equation}
\end{thm}

\begin{proof}
If $D\vartriangleleft C$ then certainly $D^*\vartriangleleft C^*$.  If $D'\vartriangleleft C'$ too and $\mathsf{s}(C)=\mathsf{r}(C')$ then \autoref{DirectedUnitCosets} and \eqref{triEq} yield $\mathsf{s}(D)=(\mathsf{s}(C)\cap N)^<=(\mathsf{r}(C')\cap N)^<=\mathsf{r}(D')$ and
\[\mathsf{r}(C\cdot C')\cap N=\mathsf{r}(C)\cap N=\mathsf{r}(D)\cap N=\mathsf{r}(D\cdot D')\cap N,\]
so $D\cdot D'\vartriangleleft C\cdot C'$.  This shows that $\vartriangleleft$ is functorial.

Now say that $\mathsf{r}(D)\vartriangleleft U$, for some unit coset $U$.  For any $d\in D$, we see that $\mathsf{r}(D)|d$ so $U|d$ and $D=(\mathsf{r}(D)d)^<\subseteq(Ud)^<$, as $\mathsf{r}(D)\subseteq U$.  By \autoref{A|b}, $(Ud)^<$ is an coset with $\mathsf{r}(Ud)=U$.  Thus $\mathsf{r}(D)\cap N=U\cap N=\mathsf{r}(Ud)\cap N$ and hence $D\vartriangleleft(Ud)^<$, by \eqref{triEq}, showing that $\vartriangleleft$ is star-surjective.  On the other hand, if we had another $C\in\mathcal{C}(S)$ with $D\vartriangleleft C$ and $\mathsf{r}(C)=U$ then $C=(\mathsf{r}(C)d)^<=(Ud)^<$, as $d\in D\subseteq C$, showing that $\vartriangleleft$ is also star-injective.

Laslty, note that if $\mathcal{D}_a\ni D\vartriangleleft C$ then certainly $C\in\mathcal{C}_a$, while conversely if $C\in\mathcal{C}_a$ then \autoref{tri} yields $D\in\mathcal{D}_a$ with $D\vartriangleleft C$.  This shows that $\vartriangleleft^{-1}[\mathcal{D}_a]=\mathcal{C}_a$ so $\vartriangleleft$ is also continuous and hence a Zakrzewski morphism satisfying \eqref{CtriD}.
\end{proof}

We call an \'etale representation $\mu:S\rightarrow\mathcal{B}(G)$ \emph{Zakrzewski-universal} if, for every \'etale representation $\theta:S\rightarrow\mathcal{B}(H)$, there exists a unique Zakrzewski morphism $\phi\subseteq G\times H$ such that $\theta=\overline\phi\circ\mu$, where $\overline\phi:\mathcal{B}(G)\rightarrow\mathcal{B}(H)$ is the preimage map.

\begin{thm}\label{ZakrzewskiUniversalEtaleRepresentation}
If $Z$ is symmetric then the directed coset representation $\mathcal{D}$ is a Zakrzewski-universal \'etale representation.
\end{thm}

\begin{proof}
As $\mathcal{C}$ is an \'etale representation, by \autoref{SymmetricEtaleRepresentation}, and $\mathcal{D}(S)$ is an ideal of $\mathcal{C}(S)$, by \autoref{DIdeal}, it follows that $\mathcal{D}$ is also an \'etale representation.

To see that $\mathcal{D}$ is Zakrzewski-universal, let us take another \'etale representation $\theta:S\rightarrow\mathcal{B}(G)$.  Let $\phi$ be the \'etale morphism in the proof of \autoref{SymmetricEtaleRepresentation} given by
\[\phi(g)=\{a\in S:g\in\theta(a)\}.\]
By \autoref{triEtaleMorphism}, $\vartriangleleft$ is also a Zakrzewski morphism and hence so is the composition
\[\psi=\ \vartriangleleft\circ\,\phi.\]
Moreover, $\overline{\psi}\circ\mathcal{D}=\overline{\phi}\circ\overline{\vartriangleleft}\circ\mathcal{D}=\overline{\phi}\circ\mathcal{C}=\theta$, by \eqref{CtriD} and the proof of \autoref{SymmetricEtaleRepresentation}.

The only thing left to prove is that $\psi$ is unique.  Accordingly, say we had another Zakrzewski morphism $\psi'\subseteq\mathcal{D}(S)\times G$ with $\overline{\psi'}\circ\mathcal{D}=\theta$, i.e. $\psi'^{-1}[\mathcal{D}_a]=\theta(a)$, for all $a\in S$.  Then we see that $g\in\theta(a)$ iff $D\ \psi'\ g$, for some $D\in\mathcal{D}_a$, and hence
\[\bigcup_{D\,\psi'\,g}D=\{a\in S:g\in\theta(a)\}=\phi(g).\]
As $\psi'$ is a functorial relation, all those $D\in\mathcal{D}(S)$ with $D\,\psi'\,g$ must have a common source.  By \autoref{UniqueCover}, it follows that
\[\{D\in\mathcal{D}(S):D\ \psi'\ g\}=\phi(g)^\vartriangleright=\{D\in\mathcal{D}(S):D\ \psi\ g\}.\]
This shows that $\psi'=\psi$, as required.
\end{proof}

As in \autoref{PullbackEtale}, we can consider the pullback bundle $\rho_\mathcal{D}^\vartriangleleft:\mathop{\vartriangleleft^{\rho_\mathcal{D}}}\widetilde{\mathcal{D}}(S)\rightarrow\mathcal{C}(S)$, which is continuously isomorphic to the coset bundle $\rho_\mathcal{C}:\widetilde{\mathcal{C}}(S)\rightarrow\mathcal{C}(S)$.

\begin{prp}\label{iotaContinuousIsomorphism}
Whenever $d\in D\vartriangleleft C$, we can define
\[\iota([d,D],C)=[d,C].\]
The resulting function is a continuous isomorphism $\iota:\mathop{\vartriangleleft^{\rho_\mathcal{D}}}\widetilde{\mathcal{D}}(S)\rightarrow\widetilde{\mathcal{C}}(S)$.
\end{prp}

\begin{proof}
If $d\in D\vartriangleleft C$ then $C^Z=D^Z$ and ${}^ZC={}^ZD$, by \eqref{CtriD=>CZ=DZ}, so $d^{\sim_C}=d^{\sim_D}$, by \eqref{ZAAZ}.  In particular, $\iota$ is well-defined.  Whenever $d\in C\in\mathcal{C}(S)$, \autoref{tri} yields a unique $D\vartriangleleft C$ with $d\in D$.  This shows that $\iota$ is a bijection, which is immediately seen to be a groupoid isomorphism.  Also, for all $s\in S$, we see that $\iota^{-1}[\widetilde{\mathcal{C}}_s]=\rho_\mathcal{D}^{\vartriangleleft-1}[\mathcal{C}_s]$ and $\iota^{-1}[\widetilde{\mathcal{C}}^s]=\mathop{\vartriangleleft^{\rho_\mathcal{D}}}\widetilde{\mathcal{D}}(S)\cap(\widetilde{\mathcal{D}}^s\times\mathcal{C}_s)$ so $\iota$ is also continuous.
\end{proof}

Now we have the following extension of \autoref{triEtaleMorphism}.  As before, let $\widetilde{\mathcal{C}}$ and $\widetilde{\mathcal{D}}$ denote the maps $a\mapsto\widetilde{a}=\widetilde{a}_\mathcal{C}$ and $a\mapsto\widetilde{a}_\mathcal{D}=\widetilde{a}|_{\mathcal{D}(S)}$ respectively.

\begin{cor}
$(\vartriangleleft,\iota)$ is a Zakrzewski-Pierce morphism from $\rho_\mathcal{D}$ to $\rho_\mathcal{C}$ such that
\begin{equation}\label{CiotatriD}
\widetilde{\mathcal{C}}=\tfrac{\iota}{\vartriangleleft}\circ\widetilde{\mathcal{D}},
\end{equation}
where $\tfrac{\iota}{\vartriangleleft}:\mathcal{S}(\rho_\mathcal{D})\rightarrow\mathcal{S}(\rho_\mathcal{C})$ is the semigroup homomorphism from \autoref{BundleMorphism->SemigroupHomo}.
\end{cor}

\begin{proof}
By \autoref{triEtaleMorphism} and \autoref{iotaContinuousIsomorphism}, $(\vartriangleleft,\iota)$ is a Zakrzewski-Pierce morphism.  For any $C\in\mathrm{dom}(\widetilde{a}_\mathcal{C})=\mathcal{C}_a$, we have a unique $D\in\mathcal{D}(S)$ with $a\in D\vartriangleleft C$, by \autoref{tri}.  Then we see that
\[\widetilde{a}_\mathcal{C}(C)=[a,C]=\iota([a,D],C)=\iota((\widetilde{a}_\mathcal{D}\circ\!\vartriangleleft)(C),C)=\tfrac{\iota}{\vartriangleleft}(\widetilde{a}_\mathcal{D})(C).\]
This shows that $\widetilde{a}_\mathcal{C}=\tfrac{\iota}{\vartriangleleft}(\widetilde{a}_\mathcal{D})$, which in turn shows that $\widetilde{\mathcal{C}}=\tfrac{\iota}{\vartriangleleft}\circ\widetilde{\mathcal{D}}$.
\end{proof}

Again, we call a bundle representation $\mu:S\rightarrow\mathcal{S}(\pi)$ \emph{Zakrzewski-universal} if, for every bundle representation $\theta:S\rightarrow\mathcal{S}(\pi')$, there exists a unique Zakrzewski-Pierce morphism $(\phi,\tau)$ from $\pi$ to $\pi'$ such that $\theta=\tfrac{\tau}{\phi}\circ\mu$.

\begin{thm}\label{ZakrzewskiUniversalBundleRepresentation}
If $Z$ is symmetric then the directed coset bundle representation $\widetilde{\mathcal{D}}$ is a Zakrzewski-universal bundle representation.
\end{thm}

\begin{proof}
As $\widetilde{\mathcal{C}}$ is a bundle representation, by \autoref{CosetBundleUniversal}, and $\mathcal{D}(S)$ is an ideal of $\mathcal{C}(S)$, by \autoref{DIdeal}, it follows that $\widetilde{\mathcal{D}}$ is also a bundle representation.

To see that $\widetilde{\mathcal{D}}$ is Zakrzewski-universal, let us take another bundle representation $\theta:S\rightarrow\mathcal{S}(\pi)$.  Let $\psi=\ \vartriangleleft\circ\,\phi$ be the Zakrzewski morphism in the proof of \autoref{ZakrzewskiUniversalEtaleRepresentation}, where $\phi$ is the \'etale morphism in the proof of \autoref{SymmetricEtaleRepresentation} given by $\phi(g)=\{a\in S:g\in\theta(a)\}$.  Further let $\tau:\phi^\pi\widetilde{\mathcal{C}}(S)\rightarrow F$ be the continuous functor in the proof of \autoref{CosetBundleUniversal} given by $\tau([a,\phi(g)],g)=\theta(a)(g)$.  We then get another continuous functor $\sigma:\psi^{\rho_\mathcal{D}}\widetilde{\mathcal{D}}(S)\rightarrow F$ defined by
\[\sigma([a,D],g)=\tau(\iota([a,D],\phi(g)),g)=\tau([a,\phi(g)],g)=\theta(a)(g)\]
whenever $a\in D\vartriangleleft\phi(g)$.  Thus $(\psi,\sigma)$ is a Zakrzewski-Pierce morphism with
\[\tfrac{\sigma}{\psi}(\widetilde{a}_\mathcal{D})(g)=\sigma((\widetilde{a}_\mathcal{D}\circ\psi)(g),g)=\sigma([a,D],g)=\theta(a)(g),\]
where $D$ is the unique directed coset with $a\in D\vartriangleleft\phi(g)$.  This shows that $\tfrac{\sigma}{\psi}(\widetilde{a}_\mathcal{D})=\theta(a)$, which in turn shows that $\theta=\tfrac{\sigma}{\psi}\circ\widetilde{\mathcal{D}}$.

For uniqueness, say we had another Zakrzewski-Pierce morphism $(\psi',\sigma')$ with $\theta=\tfrac{\sigma'}{\psi'}\circ\widetilde{\mathcal{D}}$.  By \autoref{ZakrzewskiUniversalEtaleRepresentation}, we must have $\psi'=\psi$.  For any $([a,D],g)\in\psi^{\rho_\mathcal{D}}\widetilde{\mathcal{D}}(S)$,
\[\sigma'([a,D],g)=\tfrac{\sigma'}{\psi'}(\widetilde{a}_\mathcal{D})(g)=\theta(a)(g)=\sigma([a,D],g).\]
Thus $\sigma'=\sigma$, proving uniqueness and hence universality.
\end{proof}

As in \autoref{CosetBundleFaithful}, it follows that even directed cosets and their equivalence classes determine when structured semigroups have faithful bundle representations.

\begin{cor}\label{DirectedCosetBundleFaithful}
If $Z$ is symmetric then the following are equivalent.
\begin{enumerate}
\item $(S,N,Z)$ has a faithful bundle representation.
\item The directed coset bundle representation of $(S,N,Z)$ is faithful.
\item For all distinct $a,b\in S$, we have $D\in\mathcal{D}(S)$ with $\{a,b\}\cap D\neq\emptyset$ and $a\not\sim_Db$.
\end{enumerate}
\end{cor}

\begin{proof}
Exactly like the proof of \autoref{CosetBundleFaithful}, just using \autoref{ZakrzewskiUniversalBundleRepresentation}.
\end{proof}

\bibliography{maths}{}
\bibliographystyle{alphaurl}

\end{document}